\documentclass[10pt]{amsart}
\usepackage[T1]{fontenc}
\usepackage{geometry}
\usepackage[latin1] {inputenc}
\usepackage{amsmath}
\usepackage{amsfonts, amssymb, textcomp}
\usepackage[colorlinks=flase, linkcolor=red,urlcolor=green, citecolor=blue]{hyperref}
\usepackage{subeqnarray}

\usepackage{latexsym}
\usepackage{fancyhdr}
\usepackage{longtable}
\usepackage{amsmath, amssymb}
\usepackage{graphicx}
\usepackage{enumerate}

\numberwithin{equation}{section}

\theoremstyle{plain}
\newtheorem{proposition}{Proposition}[section]
\newtheorem{theorem}[proposition]{Theorem}
\newtheorem{lemma}[proposition]{Lemma}
\newtheorem{corollary}[proposition]{Corollary}

\newtheorem{example}[proposition]{Example}

\newtheorem{remark}[proposition]{Remark}

\newcommand{\RR}{\mathbb{R}}

\newcommand{\NN}{\mathbb{N}}

\let\on=\operatorname

\newenvironment{demo}[1]{\par\smallskip\noindent{\bf #1.}}{\par\smallskip}

\makeatletter
\@namedef{subjclassname@2020}{%
  \textup{2020} Mathematics Subject Classification}
\makeatother

\newsavebox{\fmbox}
\newenvironment{fmpage}[1]
 {\begin{lrbox}{\fmbox}\begin{minipage}{#1}}
 {\end{minipage}\end{lrbox}\fbox{\usebox{\fmbox}}}

\title[On subadditivity-like conditions for associated weight functions]
{On subadditivity-like conditions for associated weight functions}

\author[G.~Schindl]{Gerhard Schindl}

\address{G.~Schindl: Fakult\"at f\"ur Mathematik, Universit\"at Wien, Oskar-Morgenstern-Platz~1, A-1090 Wien, Austria.}
\email{gerhard.schindl@univie.ac.at}

\begin{document}

\begin{abstract}
The aim of this article is to provide characterizations for subadditivity-like growth conditions for the so-called associated weight functions in terms of the defining weight sequence. Such growth requirements arise frequently in the literature and are standard when dealing with ultradifferentiable function classes defined by Braun-Meise-Taylor weight functions since they imply or even characterize important and desired consequences for the underlying function spaces, e.g. closedness under composition.
\end{abstract}

\thanks{G. Schindl is supported by FWF-Project P33417-N}
\date{\today}

\maketitle

\section{Introduction}\label{Introduction}
In the theory of ultradifferentiable function spaces there exist two classical, in general distinct (see \cite{BonetMeiseMelikhov07}), approaches in order to control the growth of the derivatives of the functions belonging to such classes: Either one uses a weight sequence $M=(M_p)_p$ or a weight function $\omega:[0,+\infty)\rightarrow[0,+\infty)$. In both settings one requires several basic growth and regularity assumptions on $M$ and $\omega$ and one distinguishes between two types, the {\itshape Roumieu-type spaces} $\mathcal{E}_{\{M\}}$ and $\mathcal{E}_{\{\omega\}}$, and the {\itshape Beurling-type spaces} $\mathcal{E}_{(M)}$ and $\mathcal{E}_{(\omega)}$. In the following we write $\mathcal{E}_{[\star]}$ if we mean either $\mathcal{E}_{\{\star\}}$ or $\mathcal{E}_{(\star)}$, but not mixing the cases.\vspace{6pt}

Subadditivity-like growth assumptions on $\omega$ play a fundamental role in the study of the classes $\mathcal{E}_{[\omega]}$.

Originally, in \cite{Bjorck66} subadditivity has been assumed for $\omega$ (condition $(\alpha)$ there) and the decay of the Fourier transform of the (compactly supported) functions under consideration has been measured in terms of $\omega$. Then, in \cite{BraunMeiseTaylor90} this approach has been rewritten and transformed into the nowadays frequently used form by controlling the growth of the derivatives. This has been possible by involving the {\itshape Young conjugate} $\varphi^{*}_{\omega}$ of $\varphi_{\omega}: t\mapsto\omega(e^t)$ (see \eqref{legendreconjugate}) and by requiring convexity for $\varphi_{\omega}$.

But, as mentioned in the introduction of \cite{BraunMeiseTaylor90}, a crucial improvement of this work was the replacement of subadditivity by the weaker assumption
$$\omega(2t)=O(\omega(t)),\;\;\;t\rightarrow+\infty.$$
In \cite{BraunMeiseTaylor90} and in many other places in the literature this growth condition is denoted by $(\alpha)$. In this paper we write \hyperlink{om1}{$(\omega_1)$} for it. This condition is by now standard. In particular, \hyperlink{om1}{$(\omega_1)$} is crucial (by absorbing exponential growth, see \eqref{newexpabsorb}) for describing $\mathcal{E}_{[\omega]}$ as a topological vector space in terms of the so-called {\itshape associated weight matrix}, a notion introduced in \cite{compositionpaper} and \cite{dissertation}.

For the sake of completeness let us mention that more recently in \cite[Cor. 2.14]{index} an equivalence between \hyperlink{om1}{$(\omega_1)$} and the positivity of a new growth index $\gamma(\omega)$ (see \eqref{newindex2}) has been given under very general basic growth assumptions on $\omega$. Note that this index, in particular its positivity, is becoming relevant when proving extension theorems in the ultraholomorphic weight function setting.\vspace{6pt}

However, also the subadditivity is still playing a crucial role. A weight function $\omega$ is equivalent to a subadditive weight function $\sigma$ (by equivalence we mean $\omega(t)=O(\sigma(t))$ and $\sigma(t)=O(\omega(t))$ as $t\rightarrow+\infty$) if and only if
\begin{equation}\label{alpha0}
\exists\;C\ge 1\;\exists\;t_0\ge 0\;\forall\;\lambda\ge 1\;\forall\;t\ge t_0:\;\;\;\omega(\lambda t)\le C\lambda\omega(t).
\end{equation}
In the literature this condition is frequently denoted by $(\alpha_0)$. It is known that $(\alpha_0)$ is characterizing for $\mathcal{E}_{[\omega]}$ having desired stability properties, e.g. closedness under composition or closedness under solving ODE's. We refer to \cite{compositionpaper}, \cite[Thm. 1, Thm. 3]{characterizationstabilitypaper} and \cite[Theorem 4.8]{almostanalytic} and the references therein (e.g. see \cite{FernandezGalbis06}). More detailed explanations are given in Section \ref{alpha01comments}.\vspace{6pt}

In \cite{PetzscheVogt} a very similar looking condition has been considered, denoted by $(\alpha_1)$ in this paper:
\begin{equation}\label{alpha1}
\sup_{\lambda\ge 1}\limsup_{t\rightarrow\infty}\frac{\omega(\lambda t)}{\lambda\omega(t)}<\infty.
\end{equation}
In \cite[Prop. 1.1]{PetzscheVogt} it has been stated that $(\alpha_0)$ is equivalent to $(\alpha_1)$ but in \cite[p. 27]{dissertation} we have shown that the weight $\widetilde{\omega}(t):=t\log(t)$ (for $t>0$ sufficiently large, see \eqref{omegacounter}) has $(\alpha_1)$ but not $(\alpha_0)$. By inspecting the arguments given in \cite[Prop. 1.1]{PetzscheVogt} we get that in this result condition $(\alpha_1)$ should be replaced by $(\alpha_0)$, see also \cite[Lemma 3.8.4]{dissertation}. Summarizing we have the implications
$$(\alpha_0)\Longrightarrow(\alpha_1)\Longrightarrow\hyperlink{om1}{(\omega_1)}.$$

A natural question is the meaning of $(\alpha_1)$ for the classes $\mathcal{E}_{[\omega]}$ and up to the best of our knowledge so far we only have the following information (which has served as the main motivation for writing this article and for which the author is grateful for the explanations given by Prof. Jos\'{e} Bonet Solves and Andreas Debrouwere in private communications):

In the main result \cite[Thm. 2.1]{BonetDomanski07} it has been shown that the space $\mathcal{E}_{\{\omega\}}(U)$, $U\subseteq\RR^d$ convex open, has the so-called {\itshape dual estimation property for small} $\theta$, which is a linear topological invariant for $(PLS)$-spaces introduced in \cite{BonetDomanski07}, provided that $\omega$ has $(\alpha_1)$. In fact, this result uses \cite[Thm. 6.2]{BonetDomanski06} where another ''$(\Omega)$-type invariant'' in the theory of $(PLS)$-spaces, namely property $(P\overline{\overline{\Omega}})$, is ensured for $\mathcal{E}_{\{\omega\}}(U)$ when $\omega$ is having $(\alpha_1)$.

However, in \cite{AdeBrouwere20} the author has been able to improve \cite[Thm. 2.1]{BonetDomanski07} by dropping the assumptions that $U$ is convex and $(\alpha_1)$.\vspace{6pt}

Now, when given a weight sequence $M$ then one can consider the {\itshape associated weight function} $\omega_M$ studied in detail in \cite[Chapitre I]{mandelbrojtbook}, see Section \ref{assofctsect}. When assuming some basic assumptions on $M$, then it is known that $\omega_M$ has already several standard assumptions in the weight function setting but \hyperlink{om1}{$(\omega_1)$} and $(\alpha_0)$ (and $(\alpha_1)$) are not guaranteed in general. The main aim of this paper is now to give characterizations of all these conditions in terms of a given $M$. In general, this seems to be a natural question in the ultradifferentiable weight sequence setting and moreover we see that these results are helpful to treat the following problems:

\begin{itemize}
\item[$(a)$] Understand in a better and more precise way the quantitative difference between the aforementioned subadditivity-like conditions for weight functions.

\item[$(b)$] Construct (counter)-examples of weight functions with prescribed subadditivity-like conditions. However, by the characterizing results Theorem \ref{alpha0theorem} and Theorem \ref{alpha0theoremcor} it turns out that the difference between $(\alpha_0)$ and $(\alpha_1)$ is subtle and in general it seems to be quite technical to provide an example such that $\omega_M$ has $(\alpha_1)$ but $(\alpha_0)$ is violated.

\item[$(c)$] In applications in the weight sequence setting the associated weight function arises naturally and the subadditivity-like growth properties from above can become crucial for the proofs.

\item[$(d)$] It is convenient to transfer known results from the weight function to the weight sequence setting, more precisely from the classes $\mathcal{E}_{[\omega_M]}$ to $\mathcal{E}_{[M]}$ (reduction arguments), which is in general possible by \cite{BonetMeiseMelikhov07} and \cite{compositionpaper}. (For example, in \cite[Sect. 5]{PetzscheVogt} this technique has been applied.) However, in order to do so the knowledge of \hyperlink{om1}{$(\omega_1)$} and $(\alpha_0)$ (and $(\alpha_1)$) for $\omega_M$ is becoming relevant.
\end{itemize}

Note that $(a)$ and $(b)$ might have applications not only within the ultradifferentiable (or ultraholomorphic) setting but also in general when working with (weight) functions.\vspace{6pt}

The paper is organized as follows: After collecting all technical properties on weight sequences, functions and matrices in Section \ref{weightcondsection}, in Section \ref{omega1charactsect} we provide a precise characterization of \hyperlink{om1}{$(\omega_1)$} in terms of $M$ (see \eqref{omega1charactequalternative} in Theorem \ref{omega1charact}) and compare it with known conditions and results. In Section \ref{alpha01sect} we give a full characterization of $(\alpha_0)$ in Theorem \ref{alpha0theorem} and in Theorem \ref{alpha0theorem1rem} when assuming in addition the technical condition {\itshape moderate growth}. We compare both statements and its consequences with known results concerning stability results for ultradifferentiable function classes. In Section \ref{alpha1section}, for the sake of completeness, we give some comments on the characterization of $(\alpha_1)$ and compare it with $(\alpha_0)$. In the final Section \ref{addcommsection} we summarize several observations concerning known conditions which are related to the index $\gamma(\omega)$ and which are even stronger than $(\alpha_0)$.\vspace{6pt}

\textbf{Acknowledgements.}

The author of this article thanks the anonymous referee for the careful reading and the valuable suggestions which have improved and clarified the presentation of this paper.

\section{Weights and conditions}\label{weightcondsection}
\subsection{General notation}
We write $\NN:=\{0,1,2,\dots\}$ and $\NN_{>0}:=\{1,2,3,\dots\}$.

\subsection{Weight sequences}\label{weightsequences}
Given a sequence $M=(M_p)_p\in\RR_{>0}^{\NN}$ we also use $m=(m_p)_p$ defined by $m_p:=\frac{M_p}{p!}$ and $\mu_p:=\frac{M_p}{M_{p-1}}$, $\mu_0:=1$, and analogously for all other arising sequences. $M$ is called {\itshape normalized} if $1=M_0\le M_1$ holds true.\vspace{6pt}

$M$ is called {\itshape log-convex} if
$$\forall\;p\in\NN_{>0}:\;M_p^2\le M_{p-1} M_{p+1},$$
equivalently if $\mu=(\mu_p)_p$ is nondecreasing. If $M$ is log-convex and normalized, then both $p\mapsto M_p$ and $p\mapsto(M_p)^{1/p}$ are nondecreasing and $(M_p)^{1/p}\le\mu_p$ for all $p\in\NN_{>0}$. Finally $M_p M_q\le M_{p+q}$ follows for all $p,q\in\NN$.\vspace{6pt}

If $m$ is log-convex, then $M$ is called {\itshape strongly log-convex}, denoted by \hypertarget{slc}{$(\text{slc})$}. (In \cite{PetzscheVogt} this requirement is denoted by $(\alpha_1)$ for $M$.) For our purpose it is convenient to consider the following set of sequences
$$\hypertarget{LCset}{\mathcal{LC}}:=\{M\in\RR_{>0}^{\NN}:\;M\;\text{is normalized, log-convex},\;\lim_{p\rightarrow+\infty}(M_p)^{1/p}=+\infty\}.$$
We see that $M\in\hyperlink{LCset}{\mathcal{LC}}$ if and only if $1=\mu_0\le\mu_1\le\dots$, $\lim_{p\rightarrow+\infty}\mu_p=+\infty$ (e.g. see \cite[p. 104]{compositionpaper}) and there is a one-to-one correspondence between $M$ and $\mu=(\mu_p)_p$ by taking $M_p:=\prod_{i=0}^p\mu_i$.\vspace{6pt}

For any $M\in\hyperlink{LCset}{\mathcal{LC}}$ and $C\in\NN_{>0}$ we get that
\begin{equation}\label{rootgeneralizedincr}
\forall\;p\in\NN:\;\;\;(M_{Cp})^{1/C}\le(M_{(C+1)p})^{1/(C+1)},
\end{equation}
because this is equivalent to having $M_{Cp}(M_{Cp})^{1/C}\le M_{Cp+p}$ and so to $\mu_1\cdots\mu_{Cp}=M_{Cp}\le(\mu_{Cp+1}\cdots\mu_{Cp+p})^C$ which is satisfied since $p\mapsto\mu_p$ is nondecreasing.\vspace{6pt}


$M$ has condition {\itshape moderate growth}, denoted by \hypertarget{mg}{$(\text{mg})$}, if
$$\exists\;C\ge 1\;\forall\;p,q\in\NN:\;M_{p+q}\le C^{p+q} M_p M_q.$$
In \cite{Komatsu73} this is denoted by $(M.2)$ and also known in the literature under the name {\itshape stability under ultradifferential operators.} This requirement does hold true simultaneously for both $M$ and $m$ by changing the constant $C$.\vspace{6pt}

We say that $M$ has condition \hypertarget{beta1}{$(\beta_1)$} (introduced in \cite{petzsche}), if
$$\exists\;Q\in\NN_{\ge 2}:\;\;\;\liminf_{p\rightarrow+\infty}\frac{\mu_{Qp}}{\mu_p}>Q.$$
In \cite{petzsche} the surjectivity of the {\itshape Borel mapping} $f\mapsto(f^{(p)}(0))_{p\in\NN}$ from $\mathcal{E}_{[M]}$ onto the corresponding weighted sequence space has been characterized in terms of this condition. Moreover, there for $M\in\hyperlink{LCset}{\mathcal{LC}}$ it has been shown that \hyperlink{beta1}{$(\beta_1)$} is equivalent to requiring \hyperlink{gamma1}{$(\gamma_1)$}, i.e.
$$\sup_{p\in\NN_{>0}}\frac{\mu_p}{p}\sum_{k\ge p}\frac{1}{\mu_k}<+\infty.$$
In the literature \hyperlink{gamma1}{$(\gamma_1)$} is also called ''strong nonquasianalyticity condition'' and in \cite{Komatsu73} it is denoted by $(M.3)$ (in fact, there $\frac{\mu_p}{p}$ is replaced by $\frac{\mu_p}{p-1}$ for $p\ge 2$ but this is equivalent to having \hyperlink{gamma1}{$(\gamma_1)$}).\vspace{6pt}


A weaker requirement than \hyperlink{beta1}{$(\beta_1)$} is
\begin{equation}\label{beta3}
\exists\;Q\in\NN_{\ge 2}:\;\;\;\liminf_{p\rightarrow+\infty}\frac{\mu_{Qp}}{\mu_p}>1,
\end{equation}
which is arising in the main characterizing results in \cite{BonetMeiseMelikhov07} and denoted by $(\beta_3)$ in \cite{dissertation}. Conditions of this type are also showing up in \cite[Thm. 3.11 $(v)$]{index}.\vspace{6pt}

Let $M,N\in\RR_{>0}^{\NN}$ be given, we write $M\hypertarget{preceq}{\preceq}N$ if $\sup_{p\in\NN_{>0}}\left(\frac{M_p}{N_p}\right)^{1/p}<+\infty$ and call $M$ and $N$ {\itshape equivalent}, denoted by $M\hypertarget{approx}{\approx}N$, if $M\hyperlink{preceq}{\preceq}N$ and $N\hyperlink{preceq}{\preceq}M$. Property \hyperlink{mg}{$(\on{mg})$} is clearly preserved under \hyperlink{approx}{$\approx$} and for \hyperlink{beta1}{$(\beta_1)$} this follows by the characterizations obtained in \cite{petzsche}.

Finally, write $M\le N$ if $M_p\le N_p$ for all $p\in\NN$.\vspace{6pt}


Concerning condition \hyperlink{mg}{$(\on{mg})$} in the literature there exists several equivalent reformulations which are useful in applications. We refer to \cite[Theorem 1]{matsumoto}, \cite[Lemma 5.3]{PetzscheVogt} and \cite[Lemma 2.2]{whitneyextensionweightmatrix} and summarize several properties needed in the forthcoming sections.

\begin{lemma}\label{mgreformulated}
Let $M\in\hyperlink{LCset}{\mathcal{LC}}$ be given, then the following are equivalent:
\begin{itemize}
\item[$(i)$] $M$ does have \hyperlink{mg}{$(\on{mg})$},

\item[$(ii)$] $M$ does satisfy $\exists\;A\ge 1\;\forall\;p\in\NN:\;\;\;M_{2p}\le A^{2p}(M_p)^2$,

\item[$(iii)$] $\mu$ does satisfy $\sup_{p\in\NN}\frac{\mu_{2p}}{\mu_p}<+\infty$,

\item[$(iv)$] $M$ does satisfy
\begin{equation}\label{mgstrange}
\exists\;A\ge 1\;\forall\;p\in\NN_{>0}:\;\;\;\mu_p\le A(M_p)^{1/p}.
\end{equation}
\end{itemize}
\end{lemma}

\eqref{mgstrange} together with log-convexity ensures that the sequences of quotients and roots are comparable up to a constant. Note that by Stirling's formula, \eqref{mgstrange} is equivalent to
\begin{equation}\label{mgstrangem}
\exists\;A\ge 1\;\forall\;p\in\NN_{>0}:\;\;\;\frac{\mu_p}{p}\le A(m_p)^{1/p}.
\end{equation}

For given $N\in\RR_{>0}^{\NN}$ and $C\in\NN_{>0}$ we set
\begin{equation}\label{sequenceLC}
L^C_p:=(N_{Cp})^{1/C},
\end{equation}
hence $L^1\equiv N$ is evident. Recall that we have put $\nu_p:=\frac{N_p}{N_{p-1}}$, $p\ge 1$, and $\nu_0:=1$ and so for the quotient sequence of $L^C$ one has
\begin{equation}\label{rootmodifiedsequence}
\lambda^C_p:=\frac{L^C_p}{L^C_{p-1}}=(\nu_{Cp-C+1}\cdots\nu_{Cp})^{1/C},\;\;\;p\ge 1,\hspace{20pt}\lambda^C_0:=1.
\end{equation}
Thus $N\in\hyperlink{LCset}{\mathcal{LC}}$ implies $L^C\in\hyperlink{LCset}{\mathcal{LC}}$ and we have the following result:

\begin{lemma}\label{LCmoderategrowth}
Let $N\in\hyperlink{LCset}{\mathcal{LC}}$ be given. Then
\begin{itemize}
\item[$(i)$] $N\le L^{C}\le L^{D}$ for any $C\le D$ and

\item[$(ii)$] $L^C\hyperlink{preceq}{\preceq}N$ holds true for some/any $C\in\NN_{\ge 2}$ if and only if $N$ has \hyperlink{mg}{$(\on{mg})$}.
\end{itemize}
Consequently, $N\hyperlink{approx}{\approx}L^C$ for some/any $C\in\NN_{\ge 2}$ if and only if $N$ has \hyperlink{mg}{$(\on{mg})$}.
\end{lemma}

\demo{Proof}
$(i)$ First, $N\le L^C$ if and only if $(N_p)^C\le N_{Cp}$ which holds true by log-convexity and normalization for $N$. The second estimate is valid by \eqref{rootgeneralizedincr}.\vspace{6pt}

$(ii)$ Suppose that $N$ has \hyperlink{mg}{$(\on{mg})$}. By Lemma \ref{mgreformulated} we get $N_{2p}\le A^{2p}(N_p)^2$ for some $A\ge 1$ and all $p\in\NN$ and by iteration
$$\forall\;x\in\NN_{>0}\;\exists\;A_x\ge 1\;\forall\;p\in\NN:\;\;\;N_{2^xp}\le A_x^{p}(N_p)^{2^x}.$$

If now $C\in\NN_{\ge 2}$ then $2^x\le C<2^{x+1}$ is valid for some $x\in\NN_{>0}$, hence $L^C_p=(N_{Cp})^{1/C}\le(N_{2^{x+1}p})^{1/(2^{x+1})}\le(A_{x+1})^{p/2^{x+1}}N_p$ since $C\mapsto(N_{Cp})^{1/C}$ is nondecreasing as seen in $(i)$. This estimates verifies $L^C\hyperlink{preceq}{\preceq}N$ and so $N\hyperlink{approx}{\approx}L^C$.\vspace{6pt}

Conversely, $L^C\hyperlink{preceq}{\preceq}N$ does imply $(N_{Cp})^{1/C}\le A^{p}N_p$ for some $A\ge 1$ and all $p\in\NN$. If $C=2$, then Lemma \ref{mgreformulated} yields directly \hyperlink{mg}{$(\on{mg})$}. If $C>2$, then $(i)$ yields $(N_{2p})^{1/2}\le(N_{Cp})^{1/C}\le A^pN_p$ and so again \hyperlink{mg}{$(\on{mg})$}.
\qed\enddemo

We close this section by recalling (a consequence of) Stirling's formula which is used frequently in the proofs of the following sections:
$$\forall\;p\in\NN_{>0}:\;\;\;\frac{p^p}{e^p}\le p!\le p^p.$$

\subsection{Associated weight function}\label{assofctsect}
Let $M\in\RR_{>0}^{\NN}$ (with $M_0=1$), then the {\itshape associated function} $\omega_M: \RR_{\ge 0}\rightarrow\RR\cup\{+\infty\}$ is defined by
\begin{equation*}\label{assofunc}
\omega_M(t):=\sup_{p\in\NN}\log\left(\frac{t^p}{M_p}\right)\;\;\;\text{for}\;t>0,\hspace{30pt}\omega_M(0):=0.
\end{equation*}
For an abstract introduction of the associated function we refer to \cite[Chapitre I]{mandelbrojtbook}, see also \cite[Definition 3.1]{Komatsu73}. If $\liminf_{p\rightarrow+\infty}(M_p)^{1/p}>0$, then $\omega_M(t)=0$ for sufficiently small $t$, since $\log\left(\frac{t^p}{M_p}\right)<0\Leftrightarrow t<(M_p)^{1/p}$ holds for all $p\in\NN_{>0}$. Moreover under this assumption $t\mapsto\omega_M(t)$ is a continuous nondecreasing function, which is convex in the variable $\log(t)$ and tends faster to infinity than any $\log(t^p)$, $p\ge 1$, as $t\rightarrow+\infty$. $\lim_{p\rightarrow+\infty}(M_p)^{1/p}=+\infty$ implies that $\omega_M(t)<+\infty$ for each finite $t$ which shall be considered as a basic assumption for defining $\omega_M$.\vspace{6pt}

If $M\in\hyperlink{LCset}{\mathcal{LC}}$, then we can compute $M$ by involving $\omega_M$ as follows, see \cite[Chapitre I, 1.4, 1.8]{mandelbrojtbook} and also \cite[Prop. 3.2]{Komatsu73}:
\begin{equation}\label{Prop32Komatsu}
M_p=\sup_{t\ge 0}\frac{t^p}{\exp(\omega_{M}(t))},\;\;\;p\in\NN.
\end{equation}
Moreover, in this case one has
\begin{equation}\label{assovanishing}
\omega_M(t)=0\hspace{20pt}\forall\;t\in[0,\mu_1],
\end{equation}
by the known integral representation formula for $\omega_M$, see \cite[1.8. III]{mandelbrojtbook} and also \cite[$(3.11)$]{Komatsu73}.


\subsection{Weight functions}\label{weightfunctions}
A function $\omega:[0,+\infty)\rightarrow[0,+\infty)$ is called a {\itshape weight function} (in the terminology of \cite[Section 2.1]{index} and \cite[Section 2.2]{sectorialextensions}), if it is continuous, nondecreasing, $\omega(0)=0$ and $\lim_{t\rightarrow+\infty}\omega(t)=+\infty$. If $\omega$ satisfies in addition $\omega(t)=0$ for all $t\in[0,1]$, then we call $\omega$ a {\itshape normalized weight function}. For convenience we will write that $\omega$ has $\hypertarget{om0}{(\omega_0)}$ if it is a normalized weight.\vspace{6pt}


Moreover we consider the following conditions, this list of properties has already been used in ~\cite{dissertation}.

\begin{itemize}
\item[\hypertarget{om1}{$(\omega_1)}$] $\omega(2t)=O(\omega(t))$ as $t\rightarrow+\infty$, i.e. $\exists\;L\ge 1\;\forall\;t\ge 0:\;\;\;\omega(2t)\le L(\omega(t)+1)$.

\item[\hypertarget{om2}{$(\omega_2)$}] $\omega(t)=O(t)$ as $t\rightarrow+\infty$.

\item[\hypertarget{om3}{$(\omega_3)$}] $\log(t)=o(\omega(t))$ as $t\rightarrow+\infty$.

\item[\hypertarget{om4}{$(\omega_4)$}] $\varphi_{\omega}:t\mapsto\omega(e^t)$ is a convex function on $\RR$.

\item[\hypertarget{om5}{$(\omega_5)$}] $\omega(t)=o(t)$ as $t\rightarrow+\infty$.

\item[\hypertarget{om6}{$(\omega_6)$}] $\exists\;H\ge 1\;\forall\;t\ge 0:\;2\omega(t)\le\omega(H t)+H$.
\end{itemize}


For convenience we define the sets
$$\hypertarget{omset0}{\mathcal{W}_0}:=\{\omega:[0,\infty)\rightarrow[0,\infty): \omega\;\text{has}\;\hyperlink{om0}{(\omega_0)},\hyperlink{om3}{(\omega_3)},\hyperlink{om4}{(\omega_4)}\},\hspace{20pt}\hypertarget{omset1}{\mathcal{W}}:=\{\omega\in\mathcal{W}_0: \omega\;\text{has}\;\hyperlink{om1}{(\omega_1)}\}.$$
For any $\omega\in\hyperlink{omset0}{\mathcal{W}_0}$ we define the {\itshape Legendre-Fenchel-Young-conjugate} of $\varphi_{\omega}$ by
\begin{equation}\label{legendreconjugate}
\varphi^{*}_{\omega}(x):=\sup\{x y-\varphi_{\omega}(y): y\ge 0\},\;\;\;x\ge 0,
\end{equation}
with the following properties, e.g. see \cite[Remark 1.3, Lemma 1.5]{BraunMeiseTaylor90}: It is convex and nondecreasing, $\varphi^{*}_{\omega}(0)=0$, $\varphi^{**}_{\omega}=\varphi_{\omega}$, $\lim_{x\rightarrow+\infty}\frac{x}{\varphi^{*}_{\omega}(x)}=0$ and finally $x\mapsto\frac{\varphi_{\omega}(x)}{x}$ and $x\mapsto\frac{\varphi^{*}_{\omega}(x)}{x}$ are nondecreasing on $[0,+\infty)$. Note that by normalization we can extend the supremum in \eqref{legendreconjugate} from $y\ge 0$ to $y\in\RR$ without changing the value of $\varphi^{*}_{\omega}(x)$ for given $x\ge 0$.

\vspace{6pt}

Let $\sigma,\tau$ be weight functions, we write $\sigma\hypertarget{ompreceq}{\preceq}\tau$ if $\tau(t)=O(\sigma(t))\;\text{as}\;t\rightarrow+\infty$
and call them equivalent, denoted by $\sigma\hypertarget{sim}{\sim}\tau$, if
$\sigma\hyperlink{ompreceq}{\preceq}\tau$ and $\tau\hyperlink{ompreceq}{\preceq}\sigma$.\vspace{6pt}

We recall the following known result, e.g. see \cite[Lemma 2.8]{testfunctioncharacterization} and \cite[Lemma 2.4]{sectorialextensions} and the references mentioned in the proofs there.

\begin{lemma}\label{assoweightomega0}
Let $M\in\hyperlink{LCset}{\mathcal{LC}}$, then $\omega_M\in\hyperlink{omset0}{\mathcal{W}_0}$ holds true. $\liminf_{p\rightarrow\infty}(m_p)^{1/p}>0$ does imply \hyperlink{om2}{$(\omega_2)$} for $\omega_M$, $\lim_{p\rightarrow\infty}(m_p)^{1/p}=+\infty$ does imply \hyperlink{om5}{$(\omega_5)$} and \hyperlink{om6}{$(\omega_6)$} for $\omega_M$ if and only if $M$ does have \hyperlink{mg}{$(\on{mg})$}.
\end{lemma}

\subsection{Weight matrices}\label{classesweightmatrices}
For the following definitions and conditions see also \cite[Section 4]{compositionpaper}.

Let $\mathcal{I}=\RR_{>0}$ denote the index set (equipped with the natural order), a {\itshape weight matrix} $\mathcal{M}$ associated with $\mathcal{I}$ is a (one parameter) family of weight sequences $\mathcal{M}:=\{M^{(x)}\in\RR_{>0}^{\NN}: x\in\mathcal{I}\}$, such that
$$\forall\;x\in\mathcal{I}:\;M^{(x)}\;\text{is normalized, nondecreasing},\;M^{(x)}\le M^{(y)}\;\text{for}\;x\le y.$$
We call a weight matrix $\mathcal{M}$ {\itshape standard log-convex,} denoted by \hypertarget{Msc}{$(\mathcal{M}_{\on{sc}})$}, if
$$\forall\;x\in\mathcal{I}:\;M^{(x)}\in\hyperlink{LCset}{\mathcal{LC}}.$$
Moreover, we put $m^{(x)}_p:=\frac{M^{(x)}_p}{p!}$ for $p\in\NN$, and $\mu^{(x)}_p:=\frac{M^{(x)}_p}{M^{(x)}_{p-1}}$ for $p\in\NN_{>0}$, $\mu^{(x)}_0:=1$.

A matrix is called {\itshape constant} if $M^{(x)}\hyperlink{approx}{\approx}M^{(y)}$ for all $x,y\in\mathcal{I}$.\vspace{6pt}

We summarize some facts which are shown in \cite[Section 5]{compositionpaper} and are needed in this work. All properties listed below are valid for $\omega\in\hyperlink{omset0}{\mathcal{W}_0}$, except \eqref{newexpabsorb} for which \hyperlink{om1}{$(\omega_1)$} is necessary (and underlines the importance of this condition in this context).

\begin{itemize}
\item[$(i)$] The idea was that to each $\omega\in\hyperlink{omset0}{\mathcal{W}_0}$ we can associate a standard log-convex weight matrix $\mathcal{M}_{\omega}:=\{W^{(l)}=(W^{(l)}_p)_{p\in\NN}: l>0\}$ by\vspace{6pt}

    \centerline{$W^{(l)}_p:=\exp\left(\frac{1}{l}\varphi^{*}_{\omega}(lp)\right)$.}\vspace{6pt}


\item[$(ii)$] $\mathcal{M}_{\omega}$ satisfies
    \begin{equation}\label{newmoderategrowth}
    \forall\;l>0\;\forall\;p,q\in\NN:\;\;\;W^{(l)}_{p+q}\le W^{(2l)}_pW^{(2l)}_q.
    \end{equation}

\item[$(iii)$] \hyperlink{om6}{$(\omega_6)$} holds if and only if some/each $W^{(l)}$ satisfies \hyperlink{mg}{$(\on{mg})$} if and only if $W^{(l)}\hyperlink{approx}{\approx}W^{(n)}$ for each $l,n>0$. Consequently \hyperlink{om6}{$(\omega_6)$} is characterizing the situation when $\mathcal{M}_{\omega}$ is constant.

\item[$(iv)$] In case $\omega$ has in addition \hyperlink{om1}{$(\omega_1)$}, then $\mathcal{M}_{\omega}$ has also
     \begin{equation}\label{newexpabsorb}
     \forall\;h\ge 1\;\exists\;A\ge 1\;\forall\;l>0\;\exists\;D\ge 1\;\forall\;p\in\NN:\;\;\;h^pW^{(l)}_p\le D W^{(Al)}_p,
     \end{equation}
and this estimate is crucial for proving $\mathcal{E}_{[\mathcal{M}_{\omega}]}=\mathcal{E}_{[\omega]}$ as locally convex vector spaces.

\item[$(v)$] We have $\omega\hyperlink{sim}{\sim}\omega_{W^{(l)}}$ for each $l>0$, more precisely
\begin{equation}\label{goodequivalenceclassic}
\forall\;l>0\,\,\exists\,D_{l}>0\;\forall\;t\ge 0:\;\;\;l\omega_{W^{(l)}}(t)\le\omega(t)\le 2l\omega_{W^{(l)}}(t)+D_{l},
\end{equation}
see \cite[Theorem 4.0.3, Lemma 5.1.3]{dissertation} and also \cite[Lemma 2.5]{sectorialextensions}.
\end{itemize}

Since for any given $M\in\hyperlink{LCset}{\mathcal{LC}}$ we have $\omega_M\in\hyperlink{omset0}{\mathcal{W}_0}$, see Lemma \ref{assoweightomega0}, it does make sense to define the matrix associated with the weight $\omega_M$ by
$$\mathcal{M}_{\omega_M}:=\{M^{(l)}: l>0\}.$$
Then we get
\begin{equation}\label{Mxonestable}
\forall\;p\in\NN:\;\;\;M_p=M^{(1)}_p,
\end{equation}
which follows by applying \eqref{Prop32Komatsu} (see also the proof of \cite[Thm. 6.4]{testfunctioncharacterization}):
 \begin{align*}
 M^{(1)}_p&:=\exp(\varphi^{*}_{\omega_{M}}(p))=\exp(\sup_{y\ge 0}\{py-\omega_{M}(e^y)\})=\sup_{y\ge 0}\exp(py-\omega_{M}(e^y))
 \\&
 =\sup_{y\ge 0}\frac{\exp(py)}{\exp(\omega_{M}(e^y))}=\sup_{t\ge 1}\frac{t^p}{\exp(\omega_{M}(t))}=\sup_{t\ge 0}\frac{t^p}{\exp(\omega_{M}(t))}=M_p.
 \end{align*}
For this recall that by normalization we have $\omega_{M}(t)=0$ for $0\le t\le 1$, see \eqref{assovanishing}.

Moreover, one has by definition
\begin{equation}\label{Mgoodtransform}
\forall\;l\in\NN_{>0}\;\forall\;p\in\NN:\;\;\;M^{(l)}_p=\exp\left(\frac{1}{l}\varphi^{*}_{\omega_{M}}(l p)\right)=(M^{(1)}_{l p})^{1/l}=(M_{l p})^{1/l},
\end{equation}
and more generally
\begin{equation}\label{Mgoodtransformgeneral}
\forall\;l\in\NN_{>0}\;\forall\;C\in\NN_{>0}\;\forall\;p\in\NN:\;\;\;M^{(lC)}_p=\exp\left(\frac{1}{lC}\varphi^{*}_{\omega_{M}}(lC p)\right)=(M^{(l)}_{Cp})^{1/C}.
\end{equation}

Finally, \eqref{goodequivalenceclassic} turns into
\begin{equation}\label{goodequivalence}
\forall\;x>0\,\,\exists\,D_{x}>0\;\forall\;t\ge 0:\;\;\;x\omega_{M^{(x)}}(t)\le\omega_{M^{(1)}}(t)=\omega_M(t)\le 2x\omega_{M^{(x)}}(t)+D_{x}.
\end{equation}

\section{Characterizing condition $(\omega_1)$}\label{omega1charactsect}
We are proving now the characterization of condition \hyperlink{om1}{$(\omega_1)$} for the associated weight function $\omega_{M}$ in terms of a condition of $M$ and so we are characterizing the case $\omega_M\in\hyperlink{omset1}{\mathcal{W}}$. Moreover, we give a comparison to related known conditions. The proof of the main statement involves some associated weight matrix technique, in fact we are using the matrix $\mathcal{M}_{\omega_M}$ from the previous section.

\begin{theorem}\label{omega1charact}
Let $M\in\hyperlink{LCset}{\mathcal{LC}}$ be given. The following conditions are equivalent:
\begin{itemize}
\item[$(i)$] $M$ satisfies
\begin{equation}\label{omega1charactequ}
\exists\;L\in\NN_{>0}\;\exists\;h>1\;\exists\;A\ge 1\;\forall\;p\in\NN:\;\;\;(M_p)^Lh^{Lp}\le AM_{Lp},
\end{equation}
i.e.
\begin{equation}\label{omega1charactequalternative}
\exists\;L\in\NN_{>0}:\;\;\;\liminf_{p\rightarrow+\infty}\frac{(M_{Lp})^{\frac{1}{Lp}}}{(M_p)^{\frac{1}{p}}}>1.
\end{equation}

\item[$(ii)$] $\omega_{M}$ satisfies
$$\exists\;h>1:\;\;\;\omega_{M}(ht)=O(\omega_{M}(t)),\;\;\;t\rightarrow+\infty,$$
resp. equivalently
$$\exists\;h>1\;\exists\;L'\ge 1\;\forall\;t\ge 0:\;\;\;\omega_{M}(ht)\le L'\omega_{M}(t)+L'.$$

\item[$(iii)$] $\omega_{M}$ satisfies \hyperlink{om1}{$(\omega_1)$}.
\end{itemize}
Please note that for \eqref{omega1charactequ} we have to choose $L\ge 2$ and the proof shows that $(i)$ and $(ii)$ hold true for the same $h>1$.
\end{theorem}

\demo{Proof}
First, $(iii)$ clearly implies $(ii)$ with $h=2$. Conversely, if $(ii)$ holds true for $h\ge 2$, then $(iii)$ immediately follows since $\omega_M$ is nondecreasing. If $1<h<2$, then we are using $n$ iterations, $n\in\NN$ chosen minimal to guarantee $h^n\ge 2$, and again the fact that $\omega_M$ is nondecreasing.\vspace{6pt}

$(ii)\Rightarrow(i)$ W.l.o.g. we can assume $L'\in\NN_{>0}$. Then, by using \eqref{Prop32Komatsu}, we get for all $p\in\NN$:
\begin{align*}
M_{L'p}&=\sup_{t\ge 0}\frac{t^{L'p}}{\exp(\omega_{M}(t))}=\sup_{s\ge 0}\frac{(hs)^{L'p}}{\exp(\omega_{M}(hs))}\ge e^{-L'}\sup_{s\ge 0}\frac{(hs)^{L'p}}{\exp(L'\omega_{M}(s))}
\\&
=e^{-L'}h^{L'p}\left(\sup_{s\ge 0}\frac{s^p}{\exp(\omega_{M}(s))}\right)^{L'}=e^{-L'}h^{L'p}(M_p)^{L'},
\end{align*}
which proves $(i)$ with $A:=e^{L'}$, $L:=L'$ and the same $h>1$.\vspace{6pt}

$(i)\Rightarrow(ii)$ Let $\mathcal{M}_{\omega_{M}}:=\{M^{(l)}: l>0\}$ be the matrix associated with $\omega_{M}$. By using the assumption and \eqref{Mgoodtransform} we have
$$\exists\;L\in\NN_{\ge 2}\;\exists\;h>1\;\exists\;A\ge 1\;\forall\;p\in\NN:\;\;\;M^{(1)}_ph^p=M_ph^p\le A^{1/L}(M_{Lp})^{1/L}=A^{1/L}M^{(L)}_p,$$
hence $\log\left(\frac{(th)^p}{M^{(L)}_p}\right)\le\log\left(\frac{t^p}{M^{(1)}_p}\right)+\frac{\log(A)}{L}$ for all $t>0$ and $p\in\NN$.

Thus we obtain by definition $\omega_{M^{(L)}}(ht)\le\omega_{M^{(1)}}(t)+\frac{\log(A)}{L}$ for all $t\ge 0$ because for $t=0$ we have $\omega_{M^{(L)}}(0)=\omega_{M^{(1)}}(0)=0$.

Finally we apply \eqref{goodequivalence} and get for all $t\ge 0$
$$\omega_{M}(ht)\le 2L\omega_{M^{(L)}}(ht)+D_L\le 2L\omega_{M^{(1)}}(t)+2\log(A)+D_L=2L\omega_{M}(t)+2\log(A)+D_L,$$
which proves $(ii)$ with $L':=\max\{2L,2\log(A)+D_L\}$ and the same $h>1$.
\qed\enddemo

\begin{remark}
In a joint paper in preparation with J. Jim\'{e}nez-Garrido and J. Sanz we prove Theorem \ref{omega1charact} for weight matrices associated with weight functions $\omega$ and even more generally for abstractly given weight matrices $\mathcal{M}$. In addition, we show that \eqref{newexpabsorb} is even equivalent to \hyperlink{om1}{$(\omega_1)$} and we also focus on consequences for the associated classes of ultradifferentiable functions but which requires some more preparation.
\end{remark}

By taking into account Theorem \ref{omega1charact} we can show now that having \hyperlink{om1}{$(\omega_1)$} for $\omega_M$  is stable under the relation \hyperlink{approx}{$\approx$}. For an alternative independent and direct proof we refer to \cite[Remark 3.3]{sectorialextensions}, see also \cite[Lemma 3.18 $(1)$]{testfunctioncharacterization}.

\begin{corollary}\label{omega1stable}
Let $M\in\hyperlink{LCset}{\mathcal{LC}}$ be given such that $\omega_M$ has \hyperlink{om1}{$(\omega_1)$}. Then \hyperlink{om1}{$(\omega_1)$} is also valid for any $\omega_N$ such that $N\in\hyperlink{LCset}{\mathcal{LC}}$ and $M\hyperlink{approx}{\approx}N$ holds true.
\end{corollary}

\demo{Proof}
$M$ has \eqref{omega1charactequalternative}, hence
$$\exists\;L\in\NN_{\ge 2}\;\exists\;\varepsilon>0\;\exists\;p_{\varepsilon}\in\NN\;\forall\;p\ge p_{\varepsilon}:\;\;\;\frac{(M_{Lp})^{1/(Lp)}}{(M_p)^{1/p}}\ge 1+\varepsilon.$$
Iteration yields then $\frac{(M_{L^xp})^{1/(L^xp)}}{(M_p)^{1/p}}\ge(1+\varepsilon)^x$ for all $x\in\NN_{>0}$ and $p\ge p_{\varepsilon}$.

The equivalence $M\hyperlink{approx}{\approx}N$ gives $C^{-1}(N_p)^{1/p}\le(M_p)^{1/p}\le C(N_p)^{1/p}$ for some $C\ge 1$ and all $p\in\NN_{>0}$. Gathering everything we have

$\exists\;L\in\NN_{>0}\;\exists\;C\ge 1\;\exists\;\varepsilon>0\;\exists\;p_{\varepsilon}\in\NN\;\forall\;p\ge p_{\varepsilon}\;\forall\;x\in\NN_{>0}:$
$$C^2\frac{(N_{L^xp})^{1/(L^xp)}}{(N_p)^{1/p}}\ge\frac{(M_{L^xp})^{1/(L^xp)}}{(M_p)^{1/p}}\ge(1+\varepsilon)^x.$$
Then choose $x$ sufficiently large to ensure $\frac{(1+\varepsilon)^x}{C^2}>1+\frac{\varepsilon}{2}$ and so \eqref{omega1charactequalternative} for $N$ is verified when taking $L':=L^x$.
\qed\enddemo

Next we give a connection to condition \eqref{beta3} arising crucially in the main comparison results in \cite{BonetMeiseMelikhov07}, in particular we refer to the discussion in \cite[Remark 15]{BonetMeiseMelikhov07}. We show that the characterization in Theorem \ref{omega1charact} is consistent with the statements in \cite{BonetMeiseMelikhov07}.

\begin{proposition}\label{beta3comparsion}
Let $M\in\hyperlink{LCset}{\mathcal{LC}}$ be given and consider the following assertions:

\begin{itemize}
\item[$(i)$]  $\mu$ satisfies
$$\exists\;Q\in\NN_{\ge 2}:\;\;\;\liminf_{p\rightarrow\infty}\frac{\mu_{Qp}}{\mu_p}>1,$$
i.e. \eqref{beta3} ($(\beta_3)$) is valid.

\item[$(ii)$] $\omega_M$ has \hyperlink{om1}{$(\omega_1)$}.

\item[$(iii)$] $M$ satisfies
$$\exists\;L\in\NN_{>0}:\;\;\;\liminf_{p\rightarrow\infty}\frac{(M_{Lp})^{\frac{1}{Lp}}}{(M_p)^{\frac{1}{p}}}>1,$$
i.e. \eqref{omega1charactequalternative}.
\end{itemize}
Then we obtain the implications $(i)\Rightarrow(ii)\Leftrightarrow(iii)$. If $M$ has in addition property \hyperlink{mg}{$(\on{mg})$}, then also $(iii)\Rightarrow(i)$ is valid and so all assertions are equivalent.
\end{proposition}

In this context, note that \hyperlink{beta1}{$(\beta_1)$} implies $(\beta_3)$. So this result applies to each strongly non-quasianalytic sequence and, in particular, the equivalence of all three assertions holds for so-called {\itshape strongly regular sequences}, see \cite[Sect. 1.1]{Thilliezdivision}. Note that for any $M\in\hyperlink{LCset}{\mathcal{LC}}$ for the expressions in \eqref{beta3} and \eqref{omega1charactequalternative} we have automatically $\liminf_{p\rightarrow+\infty}(\cdot)\ge 1$.

\demo{Proof}
$(i)\Rightarrow(ii)$ This has been shown in \cite[Lemma 12, $(2)\Rightarrow(4)$]{BonetMeiseMelikhov07}.\vspace{6pt}

$(ii)\Leftrightarrow(iii)$ This is shown in Theorem \ref{omega1charact}.\vspace{6pt}

$(iii)\Rightarrow(i)$ First, we recall that by log-convexity, normalization and \hyperlink{mg}{$(\on{mg})$} (see \eqref{mgstrange}) one has
\begin{equation}\label{beta3comparsionequ}
\exists\;A\ge 1\;\forall\;p\in\NN_{>0}:\;\;\;(M_p)^{1/p}\le\mu_p\le A(M_p)^{1/p}.
\end{equation}
Now we use similar ideas as in Corollary \ref{omega1stable} above.

By \eqref{omega1charactequalternative} there exist $L\in\NN_{\ge 2}$ and $\varepsilon>0$ such that for some $p_{\varepsilon}\in\NN$ we get $\frac{(M_{Lp})^{1/(Lp)}}{(M_p)^{1/p}}\ge 1+\varepsilon$ for all $p\ge p_{\varepsilon}$. Hence, by iterating this estimate we obtain $\frac{(M_{L^np})^{1/(L^np)}}{(M_p)^{1/p}}\ge(1+\varepsilon)^n$ for all $n\in\NN_{>0}$ and $p\ge p_{\varepsilon}$. Finally, by involving \eqref{beta3comparsionequ} we get
$$\exists\;L\in\NN_{\ge 2}\;\exists\;A\ge 1\;\exists\;\varepsilon>0\;\exists\;p_{\varepsilon}\in\NN\;\forall\;p\ge p_{\varepsilon}\;\forall\;n\in\NN_{>0}:\;\;\;A\frac{\mu_{L^np}}{\mu_p}\ge\frac{(M_{L^np})^{1/(L^np)}}{(M_p)^{1/p}}\ge(1+\varepsilon)^n,$$
and so with $n$ chosen sufficiently large (e.g. such that $\frac{(1+\varepsilon)^n}{A}\ge 1+\frac{\varepsilon}{2}$ is valid) we have verified \eqref{beta3} with $Q:=L^n$.
\qed\enddemo

Note that by using \eqref{beta3comparsionequ} we can prove $(i)\Rightarrow(iii)$ analogously as $(iii)\Rightarrow(i)$. So in general one could think, by comparing roots and quotients, that also for $(i)\Rightarrow(iii)$ assumption \hyperlink{mg}{$(\on{mg})$} is required but the arguments in the previous proof yield that it is superfluous.

In fact we construct now an example which illustrates that in general the equivalence of all assertions listed in Proposition \ref{beta3comparsion} fails, i.e. the arising $\liminf$-conditions \eqref{omega1charactequalternative} and \eqref{beta3} are falling apart and so $\omega_M$ can have \hyperlink{om1}{$(\omega_1)$} even if $(\beta_3)$ for $M$ is violated.

\begin{example}\label{beta3counterexample}
Let $(a_j)_{j\ge 1}$ be a sequence in $\NN$ such that
\begin{equation}\label{sequencea}
a_1:=1,\hspace{30pt}\frac{a_{j+1}}{a_j}\ge j+1,\hspace{30pt}\frac{a_{j+1}}{a_j-1}<\frac{a_{j+2}}{a_{j+1}-1},\;\;\;\forall\;j\ge 1.
\end{equation}

We define $M=(M_p)_{p\in\NN}$ in terms of the quotients $(\mu_p)_{p\in\NN}$ via $M_p:=\prod_{i=0}^p\mu_i$. We set $\mu_0:=1$ and
\begin{equation}\label{sequencec}
\mu_{a_1}=\mu_1:=c_1:=1,\hspace{20pt}\mu_{a_j}:=c_j:=2^{\frac{a_{j+1}}{a_j-1}}(M_{a_j-1})^{\frac{1}{a_j-1}},\;\;\;j\ge 2,
\end{equation}
and
\begin{equation}\label{sequenceb}
\mu_{p}:=\mu_{a_j}=c_j,\hspace{30pt}a_j<p\le a_{j+1}-1,\;\;\;j\ge 1.
\end{equation}

{\itshape Claim:} $M\in\hyperlink{LCset}{\mathcal{LC}}$. First, we have $1=M_0=M_1$ (normalization) by $\mu_0=\mu_1=1$. In order to verify log-convexity, we have to show that $j\mapsto c_j$ is nondecreasing. Clearly, $\mu_p=\mu_{p+1}$ if $a_j\le p<p+1\le a_{j+1}-1$, $j\ge 1$, is valid and to conclude it is sufficient to verify $(M_{a_j-1})^{\frac{1}{a_j-1}}\le(M_{a_{j+1}-1})^{\frac{1}{a_{j+1}-1}}$.

By \eqref{sequencec} and \eqref{sequenceb} one has $1=\mu_1=\dots=\mu_{a_2-1}$ and so $1=M_{a_2-1}=(M_{a_2-1})^{1/(a_2-1)}$ which yields $c_2=2^{\frac{a_{3}}{a_2-1}}>1=c_1$. In general, note that $(M_p)^{1/p}\le(M_{p+1})^{1/(p+1)}$ if and only if $(M_p)^{1/p}\le\mu_{p+1}$. Thus $(M_{a_2-1})^{1/(a_2-1)}\le(M_{a_2})^{1/a_2}\le(M_{a_3-1})^{1/(a_3-1)}$ follows immediately and so $c_3>c_2$, too. Then proceed by induction. Note that this also shows $\lim_{j\rightarrow+\infty}c_j=\lim_{p\rightarrow+\infty}\mu_p=+\infty$ which proves the claim (see \cite[p. 104]{compositionpaper}).\vspace{6pt}

{\itshape Claim:} $M$ does not have $(\beta_3)$. Take $Q\in\NN_{>0}$, $Q\ge 2$, arbitrary but from now on fixed. For all $j$ with $Q\le j$ we get $a_jQ\le a_jj\le a_{j+1}-1$ because $a_jj+1\le a_j(j+1)\le a_{j+1}$ (see \eqref{sequencea}) and so by taking into account \eqref{sequenceb} we get for all $j\ge Q$
$$\frac{\mu_{Qa_j}}{\mu_{a_j}}=\frac{c_j}{c_j}=1.$$
\vspace{6pt}

{\itshape Claim:} $M$ does have \eqref{omega1charactequalternative} with $L=2$.

{\itshape Case I:} Let $a_j\le p<2p\le a_{j+1}-1$ for some $j\ge 2$. (If $j=1$, then $(M_p)^{1/p}=(M_{2p})^{1/(2p)}=1$ and so this situation has to be excluded.) Then by taking into account $M_p=\mu_1\cdots\mu_p$ we get
\begin{align*}
\frac{(M_{2p})^{1/(2p)}}{(M_p)^{1/p}}&=\frac{(M_{a_j})^{1/(2p)}(\mu_{a_j+1}\cdots\mu_{2p})^{1/(2p)}}{(M_{a_j})^{1/p}(\mu_{a_j+1}\cdots\mu_p)^{1/p}}=(M_{a_j})^{\frac{p-2p}{2p^2}}(c_j)^{(2p-a_j)/(2p)}(c_j)^{(a_j-p)/p}
\\&
=(M_{a_j})^{\frac{-1}{2p}}(c_j)^{\frac{a_j}{2p}}.
\end{align*}
Hence $\frac{(M_{2p})^{1/(2p)}}{(M_p)^{1/p}}\ge 2$ is equivalent to $(c_j)^{a_j}\ge 4^pM_{a_j}=4^p\mu_1\cdots\mu_{a_j}$ and so to
$$c_j\ge 4^{\frac{p}{a_j-1}}(\mu_1\cdots\mu_{a_j-1})^{\frac{1}{a_j-1}}=4^{\frac{p}{a_j-1}}(M_{a_j-1})^{\frac{1}{a_j-1}}.$$
Since $p\le\frac{a_{j+1}-1}{2}$ we have $4^{\frac{p}{a_j-1}}\le 2^{\frac{a_{j+1}}{a_j-1}}$ and so the desired estimate is valid by \eqref{sequencec} and \eqref{sequenceb}.\vspace{6pt}

{\itshape Case II:} If $a_j\le p\le a_{j+1}-1<2p$ for some $j\ge 1$, then $2p\le 2a_{j+1}-2\le a_{j+2}-1$ because $2a_{j+1}-1\le a_{j+2}$ by \eqref{sequencea}. In this situation, we follow the computations from case I and by $c_{j+1}>c_j$ we are able to conclude.
\end{example}

In the next result we see that \eqref{omega1charactequalternative} is weaker than the following condition: We say that $M$ is {\itshape root almost increasing}, if
\begin{equation}\label{rootalmost}
\exists\;A\ge 1\;\forall\;1\le p\le q:\;\;\;(m_p)^{1/p}\le A(m_q)^{1/q},
\end{equation}
which means that the sequence $((m_p)^{1/p})_{p\ge 1}$ is {\itshape almost increasing.} This requirement is becoming crucial when proving desired and important stability properties for the ultradifferentiable classes $\mathcal{E}_{[M]}$, e.g. closedness under composition, see \cite{compositionpaper} and \cite{characterizationstabilitypaper} and the references therein.

\begin{lemma}\label{rootalmostimplication}
Let $M\in\hyperlink{LCset}{\mathcal{LC}}$ be given. Then \eqref{rootalmost} does imply
\begin{itemize}
\item[$(i)$] \eqref{omega1charactequalternative}, i.e. \hyperlink{om1}{$(\omega_1)$} for $\omega_M$,

\item[$(ii)$] $\liminf_{p\rightarrow\infty}(m_p)^{1/p}>0$ which is yielding \hyperlink{om2}{$(\omega_2)$} for $\omega_M$ (see Lemma \ref{assoweightomega0}).
\end{itemize}
\end{lemma}

\demo{Proof}
$(i)$ Let $p\ge 1$ and apply \eqref{rootalmost} to $q=Bp$ and $B\ge 1$ is chosen such that $B\ge Ae(1+\varepsilon)$ holds true, with $A$ the constant arising in \eqref{rootalmost} and $\varepsilon>0$ arbitrary but fixed. W.l.o.g. we assume $B\in\NN_{>0}$ and then, multiplying \eqref{rootalmost} by $p!^{1/p}$, we obtain by Stirling's formula
$$(M_p)^{1/p}\le A\frac{p!^{1/p}}{(Bp)!^{1/(Bp)}}(M_{Bp})^{1/(Bp)}\le A\frac{ep}{Bp}(M_{Bp})^{1/(Bp)},$$
and so
$$\exists\;B\in\NN_{>0}\;\forall\;p\in\NN_{>0}:\;\;\;\frac{(M_{Bp})^{1/(Bp)}}{(M_p)^{1/p}}\ge\frac{B}{eA}\ge 1+\varepsilon,$$
which proves \eqref{omega1charactequalternative} with $L=B$.\vspace{6pt}

$(ii)$ \eqref{rootalmost} implies $(m_p)^{1/p}\ge A^{-1}m_1=A^{-1}M_1\ge A^{-1}$ for all $p\ge 1$ which yields the assertion by $\liminf_{p\rightarrow\infty}(m_p)^{1/p}\ge A^{-1}>0$.
\qed\enddemo

\section{The condition $(\alpha_0)$}\label{alpha01sect}
\subsection{General comments}\label{alpha01comments}
As it has already been mentioned in the introduction we have that $(\alpha_0)$ implies \hyperlink{om1}{$(\omega_1)$} (by choosing $\lambda=2$) and $(\alpha_0)$ is crucial for the characterization of several important stability properties for the ultradifferentiable classes $\mathcal{E}_{[\omega]}$. Furthermore, $(\alpha_0)$ holds if and only if the weight is equivalent to its {\itshape least concave majorant,} see \cite[Lemma 1]{peetre} and \cite[Prop. 2.23]{index}. Note that for this equivalence it is sufficient to consider weight functions as defined in Section \ref{weightfunctions} but for any (weight) function concavity and normalization can never be satisfied at the same time.\vspace{6pt}

A similar technical problem arises for $(\alpha_0)$: When $\omega(t)=0$ for all $t\in[0,t_0]$, $t_0>0$, then $C\lambda\omega(t)=0$ for any $C,\lambda\ge 1$ and $0<t\le t_0$ but $\omega(\lambda t)\rightarrow+\infty$ as $\lambda\rightarrow+\infty$. In particular, this happens for $\omega\equiv\omega_M$, $M\in\hyperlink{LCset}{\mathcal{LC}}$, see \eqref{assovanishing}. From a technical point of view one could think to switch to an equivalent but non-normalized weight. However, for the proofs, in order to apply \eqref{Prop32Komatsu}, we require the function $\omega_M$ itself and so we have to take into account the discussion just before \cite[Prop. 2.23]{index} for allowing $0$ for small values $t\ge 0$.\vspace{6pt}

Thus, for the convenience of the reader and to make everything precise we repeat the proof of \cite[Prop. 1.1]{PetzscheVogt} by using \cite[Prop. 2.23]{index} in full details and get the correct condition for all $t\ge 0$:

\begin{lemma}\label{alpha0equivalence}
Let $\omega:[0,+\infty)\rightarrow[0,+\infty)$ be nondecreasing such that $\omega(0)\ge 0$. Then the following are equivalent:
\begin{itemize}
\item[$(i)$] $\omega$ is equivalent to a subadditive function $\sigma$,

\item[$(ii)$] $\omega$ satisfies
\begin{equation}\label{alpha0rewritten}
\exists\;C\ge 1\;\exists\;D\ge 1\;\forall\;\lambda\ge 1\;\forall\;t\ge 0:\;\;\;\omega(\lambda t)\le C\lambda\omega(t)+D\lambda,
\end{equation}
\item[$(iii)$] $\omega$ is equivalent to its least concave majorant $F_{\omega}$.
\end{itemize}
\end{lemma}

In fact, \eqref{alpha0rewritten} is extending $(\alpha_0)$ to all $t\ge 0$ and is becoming important in the proofs below.

\demo{Proof}
$(i)\Rightarrow(ii)$ By equivalence we have $\omega(t)\le A\sigma(t)+A$ and $\sigma(t)\le A\omega(t)+A$ for some $A\ge 1$ and all $t\ge 0$. Take $\lambda\in\NN_{>0}$, then
$$\omega(\lambda t)\le A\sigma(\lambda t)+A\le A\lambda\sigma(t)+A\le\lambda A^2\omega(t)+A^2\lambda+A,$$
proving \eqref{alpha0rewritten} with $C=A^2$ and $D=A^2+A$ both not depending on given $\lambda$. When $\lambda\ge 1$ arbitrary, then $\mu\le\lambda<\mu+1$ for some $\mu\in\NN_{>0}$ and so
$$\omega(\lambda t)\le\omega((\mu+1)t)\le(\mu+1)A^2\omega(t)+A^2(\mu+1)+A\le2A^2\mu\omega(t)+2\mu A^2+A\le 2A^2\lambda\omega(t)+2A^2\lambda+A,$$
showing \eqref{alpha0rewritten} with $C=2A^2$ and $D=2A^2+A$.\vspace{6pt}

$(ii)\Rightarrow(iii)$ We check that \eqref{alpha0rewritten} yields \cite[$(12)$]{index}: The first condition there is clear for $C=1$ since $\omega$ is nondecreasing. For the second, we fix $s\ge 0$ and set $t:=\lambda s$ for $\lambda>1$, then \eqref{alpha0rewritten} implies $\omega(t)\le C\frac{t}{s}\omega(s)+D\frac{t}{s}$ which gives the second part of \cite[$(12)$]{index}. Thus \cite[Prop. 2.23 $(13)$]{index} (one should read $A\le 1$ there) yields the conclusion.\vspace{6pt}

$(iii)\Rightarrow(i)$ Since $F_{\omega}(0)=\omega(0)\ge 0$ (for the first equality see again \cite[Prop. 2.23]{index}), we have that $F_{\omega}$ is subadditive (e.g. see \cite[Lemma 3.8.1 $(1)$]{dissertation} for this).
\qed\enddemo

Next we are recalling for an abstractly given weight matrix $\mathcal{M}:=\{M^{(l)}: l\in\mathcal{I}\}$ the following two crucial conditions from \cite{compositionpaper} and \cite{characterizationstabilitypaper} which are generalizing \eqref{rootalmost} to a mixed setting:

\hypertarget{R-rai}{$(\mathcal{M}_{\{\text{rai}\}})$} \hskip1cm $\forall\;x\in\mathcal{I}\;\exists\;C>0\;\exists\;y\in\mathcal{I}\;\forall\;1\le p\le q:\;\;\;(m^{(x)}_p)^{1/p}\le C(m^{(y)}_q)^{1/q}$,\par\vskip.3cm

\hypertarget{B-rai}{$(\mathcal{M}_{(\text{rai})})$} \hskip1cm $\forall\;x\in\mathcal{I}\;\exists\;C>0\;\exists\;y\in\mathcal{I}\;\forall\;1\le p\le q:\;\;\;(m^{(y)}_p)^{1/p}\le C(m^{(x)}_q)^{1/q}$.\par\vskip.3cm

Here the abbreviation ''rai'' is standing for {\itshape root almost increasing}, the first requirement is the ''Roumieu-type condition'' and the second one the ''Beurling-type''. These conditions are becoming crucial for the characterization of important stability properties for ultradifferentiable classes defined in terms of weight matrices $\mathcal{E}_{\{\mathcal{M}\}}$ resp. $\mathcal{E}_{(\mathcal{M})}$ (of the particular type), see \cite[Thm. 4.9, Thm. 4.11]{compositionpaper} and \cite[Thm. 5, Thm. 6]{characterizationstabilitypaper}.

For our considerations only the Roumieu-type condition is becoming relevant. In the next result we prove that \hyperlink{R-rai}{$(\mathcal{M}_{\{\on{rai}\}})$} for the matrix $\mathcal{M}_{\omega_M}$ is equivalent to a ''truncated'' version of \eqref{rootalmost}. Note that in \hyperlink{R-rai}{$(\mathcal{M}_{\{\on{rai}\}})$} by the order $M^{(x)}\le M^{(y)}\Leftrightarrow m^{(x)}\le m^{(y)}$ for $x\le y$ we can restrict w.l.o.g. to $x,y\in\NN_{>0}$.

\begin{lemma}\label{truncatedrai}
Let $M\in\hyperlink{LCset}{\mathcal{LC}}$ be given and $\mathcal{M}_{\omega_M}:=\{M^{(l)}: l>0\}$ be the matrix associated with the weight $\omega_M$. Then the following are equivalent:
\begin{itemize}
\item[$(i)$] $\mathcal{M}_{\omega_M}$ satisfies \hyperlink{R-rai}{$(\mathcal{M}_{\{\on{rai}\}})$},

\item[$(ii)$] $m$ satisfies
\begin{equation}\label{raimixed}
\exists\;A\ge 1\;\exists\;C\in\NN_{>0}\;\forall\;p\ge 1\;\forall\;q\ge Cp:\;\;\;(m_p)^{1/p}\le A(m_{q})^{1/q}.
\end{equation}
\end{itemize}
\end{lemma}

\demo{Proof}
$(i)\Rightarrow(ii)$ By \eqref{Mxonestable} we have $m^{(1)}\equiv m$. Then we apply \hyperlink{R-rai}{$(\mathcal{M}_{\{\on{rai}\}})$} to $x=1$ and by using Stirling's formula and \eqref{Mgoodtransform} we get for some $y\in\NN_{>0}$ and all $1\le p\le q$ that
\begin{align*}
(m_p)^{1/p}&=(m^{(1)}_p)^{1/p}\le C(m^{(y)}_q)^{1/q}=C\frac{1}{q!^{1/q}}(M^{(y)}_q)^{1/q}=C\frac{1}{q!^{1/q}}(M^{(1)}_{yq})^{1/(yq)}
\\&
=C\frac{(yq)!^{1/(yq)}}{q!^{1/q}}(m^{(1)}_{yq})^{1/(yq)}\le Cey(m^{(1)}_{yq})^{1/(yq)}=Cey(m_{yq})^{1/(yq)},
\end{align*}
proving \eqref{raimixed} for $A:=Cey$ and for all $p\ge 1$ and all $q\ge yp$ such that $q=yq_1$ for some $q_1\ge p\ge 1$. In order to verify \eqref{raimixed} with $C=y$ we take $q\ge 1$ such that $yq_1<q<y(q_1+1)$ for some $q_1\ge p$, and follow a trick established in \cite{Siddiqi}, see also \cite[Theorem 4.9, $(2)\Rightarrow(3)$]{compositionpaper} (''Siddiqi trick''): For the moment set $A_1:=Cey$ and since $p\mapsto(M_p)^{1/p}$ is nondecreasing (by log-convexity and normalization) we estimate as follows:
\begin{align*}
(m_q)^{1/q}&\ge(m_{yq_1})^{1/(yq_1)}\frac{(yq_1)!^{1/(yq_1)}}{q!^{1/q}}\ge\frac{1}{A_1}(m_p)^{1/p}\frac{(yq_1)!^{1/(yq_1)}}{q!^{1/q}}\ge\frac{1}{eA_1}(m_p)^{1/p}\frac{yq_1}{q}
\\&
\ge\frac{1}{eA_1}(m_p)^{1/p}\frac{yq_1}{y(q_1+1)}=\frac{1}{eA_1}\frac{q_1}{q_1+1}(m_p)^{1/p}\ge\frac{1}{2eA_1}(m_p)^{1/p}.
\end{align*}
Consequently, \eqref{raimixed} is verified for $A:=2eA_1=2e^2Cy$ and $C:=y$.\vspace{6pt}

$(ii)\Rightarrow(i)$ First note that \eqref{Mgoodtransformgeneral} can be transferred to the sequences $m^{(l)}$ as follows (by applying Stirling's formula):
$$e^{-p}p^pm^{(lC)}_p\le p!m^{(lC)}_p=M^{(lC)}_p=(M^{(l)}_{Cp})^{1/C}=(m^{(l)}_{Cp})^{1/C}(Cp)!^{1/C}\le(m^{(l)}_{Cp})^{1/C}(Cp)^{p},$$
hence
\begin{equation}\label{mgoodtransformgeneral}
\forall\;l\in\NN_{>0}\;\forall\;C\in\NN_{>0}\;\forall\;p\in\NN_{>0}:\;\;\;(m^{(lC)}_p)^{1/p}\le eC(m^{(l)}_{Cp})^{1/(Cp)},
\end{equation}
and
$$e^{-{p}}(Cp)^{p}(m^{(l)}_{Cp})^{1/C}\le (Cp)!^{1/C}(m^{(l)}_{Cp})^{1/C}=(M^{(l)}_{Cp})^{1/C}=M^{(lC)}_{p}=p!m^{(lC)}_{p}\le p^pm^{(lC)}_{p},$$
which implies
\begin{equation}\label{mgoodtransformgeneral1}
\forall\;l\in\NN_{>0}\;\forall\;C_1\in\NN_{>0}\;\forall\;p\in\NN_{>0}:\;\;\;(m^{(l)}_{C_1p})^{1/(C_1p)}\le\frac{e}{C_1}(m^{(lC_1)}_{p})^{1/p}.
\end{equation}
Take $D\in\NN_{>0}$, arbitrary but from now on fixed. We use assumption \eqref{raimixed} with arising constants $A$ and $C$ and \eqref{mgoodtransformgeneral} and \eqref{mgoodtransformgeneral1}; the first applied to $l=1$ and $C=D$ and the second one applied to $l=1$ and $C_1=CD$. We estimate as follows for all $q\ge p\ge 1$:
\begin{align*}
(m^{(D)}_p)^{1/p}&\le eD(m^{(1)}_{Dp})^{1/(Dp)}=eD(m_{Dp})^{1/(Dp)}\le eDA(m_{CDq})^{1/(CDq)}=eDA(m^{(1)}_{CDq})^{1/(CDq)}
\\&
\le eDA e\frac{1}{CD}(m^{(CD)}_q)^{1/q}=e^2\frac{A}{C}(m^{(CD)}_q)^{1/q}.
\end{align*}
Thus we have verified \hyperlink{R-rai}{$(\mathcal{M}_{\{\on{rai}\}})$} for the matrix $\mathcal{M}_{\omega_M}$.
\qed\enddemo

We give some more information concerning condition \eqref{raimixed} by involving the auxiliary sequence $L^C$ from \eqref{sequenceLC}.

\begin{lemma}\label{truncatedraiviaLC}
Let $M\in\hyperlink{LCset}{\mathcal{LC}}$ be given. Then the following are equivalent:
\begin{itemize}
\item[$(i)$] $m$ satisfies \eqref{raimixed},

\item[$(ii)$] one has
\begin{equation}\label{raimixedcorrect}
    \exists\;B\ge 1\;\;\exists\;C\in\NN_{>0}\;\forall\;1\le p\le q:\;\;\;(m_p)^{1/p}\le B(l^C_q)^{1/q},
    \end{equation}
    with recalling $l^C_p:=\frac{L^C_p}{p!}=\frac{(M_{Cp})^{1/C}}{p!}=\frac{(Cp)!^{1/C}}{p!}(m_{Cp})^{1/C}$.
\end{itemize}
\end{lemma}

\eqref{raimixedcorrect} can be viewed as mixed \eqref{rootalmost} (resp. \hyperlink{R-rai}{$(\mathcal{M}_{\{\on{rai}\}})$} and \hyperlink{B-rai}{$(\mathcal{M}_{(\on{rai})})$}) between the two fixed weights $M$ and $L^C$.

\demo{Proof}
First, \eqref{rootmodifiedsequence} yields that $M\in\hyperlink{LCset}{\mathcal{LC}}$ implies $L^C\in\hyperlink{LCset}{\mathcal{LC}}$. Second, $l^C_p=\frac{(Cp)!^{1/C}}{p!}(m_{Cp})^{1/C}$ and so
    $$\frac{(Cp)^p}{e^pp^p}(m_{Cp})^{1/C}\le l^C_p\le\frac{e^p(Cp)^p}{p^p}(m_{Cp})^{1/C},$$
    hence
    \begin{equation}\label{raimixedcorrect1}
    \forall\;p\ge 1:\;\;\;\frac{C}{e}(m_{Cp})^{1/(Cp)}\le(l^C_p)^{1/p}\le e C(m_{Cp})^{1/(Cp)}.
    \end{equation}

$(i)\Rightarrow(ii)$ Given $p\ge 1$, then \eqref{raimixed}, when being satisfied for constants $A$ and $C$ and applied to $q=Cq'$, $q'\ge p$ arbitrary, together with the first estimate in \eqref{raimixedcorrect1} do imply \eqref{raimixedcorrect} with the same $C$ and with $B:=\frac{Ae}{C}$ (note that w.l.o.g. $A\ge C$).\vspace{6pt}

$(ii)\Rightarrow(i)$ \eqref{raimixedcorrect} and the second estimate in \eqref{raimixedcorrect1} do imply \eqref{raimixed} with $A:=BeC$ and the same $C$ for all $q=Cp'$, $p'\ge p$. If $Cp'\le q<C(p'+1)$, then apply the ''Siddiqi trick'' from Lemma \ref{truncatedrai} in order to conclude.
\qed\enddemo

By combining Lemmas \ref{LCmoderategrowth}, \ref{truncatedrai} and \ref{truncatedraiviaLC} we immediately get the following characterization.

\begin{corollary}\label{lccorollary}
Let $M\in\hyperlink{LCset}{\mathcal{LC}}$ be given and satisfying \hyperlink{mg}{$(\on{mg})$}. Then the following are equivalent:
\begin{itemize}
\item[$(i)$] $m$ satisfies \eqref{raimixed},

\item[$(ii)$] $m$ satisfies \eqref{rootalmost} (i.e. $M$ is root almost increasing).

\item[$(iii)$] $\mathcal{M}_{\omega_M}$ satisfies \hyperlink{R-rai}{$(\mathcal{M}_{\{\on{rai}\}})$},
\end{itemize}
\end{corollary}

\demo{Proof}
By Lemma \ref{LCmoderategrowth} we have $M\hyperlink{approx}{\approx}L^C$ for some/any $C\in\NN_{\ge 2}$ if and only if $M$ has \hyperlink{mg}{$(\on{mg})$} and this implies $(i)\Leftrightarrow(ii)$.
\qed\enddemo

\subsection{Main result}\label{mainresult}
In the main result we are able to give a complete characterization of $(\alpha_0)$ for $\omega_M$ in terms of given $M$.

\begin{theorem}\label{alpha0theorem}
Let $M\in\hyperlink{LCset}{\mathcal{LC}}$ be given and $\omega_M$ be the corresponding associated weight function, recall that $m_p=M_p/p!$. Then the following are equivalent:
\begin{itemize}
\item[$(i)$] $\omega_M$ satisfies
\begin{equation}\label{alpha0mixed}
\exists\;D\ge 1\;\exists\;C\ge 1\;\forall\;\lambda\ge 1\;\forall\;t\ge 0:\;\;\;\omega_M(\lambda t)\le C\lambda\omega_M(t)+D\lambda,
\end{equation}
i.e. $(\alpha_0)$ (\eqref{alpha0rewritten}).

\item[$(ii)$] $M$ satisfies
\begin{equation}\label{raimixedM}
\exists\;B\ge 1\;\exists\;C\in\NN_{>0}\;\forall\;\lambda\in\NN_{>0}\;\forall\;p\ge 1:\;\;\;(M_p)^{1/p}\lambda\le B(M_{\lambda Cp})^{1/(\lambda Cp)}.
\end{equation}

\item[$(iii)$] $m$ satisfies
$$\exists\;A\ge 1\;\exists\;C\in\NN_{>0}\;\forall\;p\ge 1\;\forall\;q\ge Cp:\;\;\;(m_p)^{1/p}\le A(m_{q})^{1/q},$$
i.e. \eqref{raimixed}.

\item[$(iv)$] The matrix $\mathcal{M}_{\omega_M}$ satisfies \hyperlink{R-rai}{$(\mathcal{M}_{\{\on{rai}\}})$}.
\end{itemize}
\end{theorem}

In particular, this result applies for any $M\in\hyperlink{LCset}{\mathcal{LC}}$ satisfying \eqref{rootalmost} because then \eqref{raimixed} holds true with $C=1$. Especially, if $M$ is satisfying \hyperlink{slc}{$(\on{slc})$}, then $p\mapsto(m_p)^{1/p}$ is nondecreasing and so we get assertion $(iii)$ with $A=C=1$. Recall that under this assumption on $M$ assertion $(i)$ has been shown in \cite[Lemma 5.5]{PetzscheVogt}.

\demo{Proof}
$(i)\Rightarrow(ii)$ W.l.o.g. we can assume $C\in\NN_{>0}$ and apply this estimate to $\lambda\in\NN_{>0}$ arbitrary but from now on fixed. Then, by using \eqref{Prop32Komatsu}, we get for all $p\in\NN$:
\begin{align*}
(M_{C\lambda p})^{1/C}&=\sup_{t\ge 0}\left(\frac{t^{C\lambda p}}{\exp(\omega_{M}(t))}\right)^{1/C}=\sup_{t\ge 0}\frac{t^{\lambda p}}{\exp(\frac{1}{C}\omega_{M}(t))}=\sup_{s\ge 0}\frac{(\lambda s)^{\lambda p}}{\exp(\frac{1}{C}\omega_{M}(\lambda s))}
\\&
\ge e^{-D\lambda/C}\sup_{s\ge 0}\frac{(\lambda s)^{\lambda p}}{\exp(\lambda\omega_{M}(s))}=e^{-D\lambda/C}\lambda^{\lambda p}\sup_{s\ge 0}\left(\frac{s^p}{\exp(\omega_{M}(s))}\right)^{\lambda}
\\&
=e^{-D\lambda/C}\lambda^{\lambda p}(M_p)^{\lambda}.
\end{align*}
So far we have shown
$$\exists\;D\ge 1\;\exists\;C\in\NN_{>0}\;\forall\;\lambda\in\NN_{>0}\;\forall\;p\in\NN:\;\;\;(M_p)^{\lambda C}\lambda^{\lambda Cp}\le e^{D\lambda}M_{\lambda Cp},$$
and so
$$\exists\;D\ge 1\;\exists\;C\in\NN_{>0}\;\forall\;\lambda\in\NN_{>0}\;\forall\;p\in\NN_{>0}:\;\;\;(M_p)^{1/p}\lambda\le e^{D/(Cp)}(M_{\lambda Cp})^{1/(\lambda Cp)},$$
proving \eqref{raimixedM} with $B:=e^D$ and $C=C$.\vspace{6pt}

$(ii)\Rightarrow(iii)$ First, dividing \eqref{raimixedM} by $p!^{1/p}$ and using Stirling's formula we arrive at
$$(m_p)^{1/p}\lambda\le B\frac{(C\lambda p)!^{1/(C\lambda p)}}{p!^{1/p}}(m_{C\lambda p})^{1/(C\lambda p)}\le B\frac{e\lambda Cp}{p}(m_{\lambda Cp})^{1/(\lambda Cp)},$$
hence
$$\exists\;B\ge 1\;\exists\;C\in\NN_{>0}\;\forall\;\lambda\in\NN_{>0}\;\forall\;p\in\NN_{>0}:\;\;\;(m_p)^{1/p}\le BeC(m_{\lambda Cp})^{1/(\lambda Cp)},$$
i.e. \eqref{raimixed} is verified with $A=BeC$ and $C=C$ for all $p,q\in\NN_{>0}$ such that $q\ge Cp$ and $q=\lambda Cp$ for some $\lambda\in\NN_{>0}$. For the remaining cases we use again the ''Siddiqi trick'' as in Lemma \ref{truncatedrai}: We set $A_1:=BeC$ and for arbitrary $q\ge Cp$ we choose $\lambda\in\NN_{>0}$ such that $\lambda Cp\le q<(\lambda+1)Cp$. By Stirling's formula and since $p\mapsto(M_p)^{1/p}$ is nondecreasing we can estimate as follows:
\begin{align*}
(m_q)^{1/q}&\ge(m_{\lambda Cp})^{1/(\lambda Cp)}\frac{(\lambda Cp)!^{1/(\lambda Cp)}}{q!^{1/q}}\ge\frac{1}{A_1}(m_p)^{1/p}\frac{(\lambda Cp)!^{1/(\lambda Cp)}}{q!^{1/q}}\ge\frac{1}{eA_1}(m_p)^{1/p}\frac{\lambda Cp}{q}
\\&
\ge\frac{1}{eA_1}(m_p)^{1/p}\frac{\lambda Cp}{(\lambda+1)Cp}=\frac{1}{eA_1}\frac{\lambda}{\lambda+1}(m_p)^{1/p}\ge\frac{1}{2eA_1}(m_p)^{1/p}.
\end{align*}
Summarizing, \eqref{raimixed} is verified for $A:=2eA_1=2e^2BC$ and $C=C$.\vspace{6pt}

$(iii)\Leftrightarrow(iv)$ This is Lemma \ref{truncatedrai}.\vspace{6pt}

$(iii)\Rightarrow(i)$ First, let us see that \eqref{raimixed} yields both \hyperlink{om1}{$(\omega_1)$} and \hyperlink{om2}{$(\omega_2)$} for $\omega_M$: In \eqref{raimixed} we take $q=BCp$, $p\ge 1$, with chosen $B\in\NN_{>0}$ such that $B>\frac{A}{C}e(1+\varepsilon)$ holds true (with $A$ and $C$ denoting the constants in \eqref{raimixed} and $\varepsilon>0$ arbitrary, but fixed). Then following the proof of Lemma \ref{rootalmostimplication} we get that \eqref{omega1charactequalternative} is verified with $L:=BC$ and so \hyperlink{om1}{$(\omega_1)$} for $\omega_M$ is follows.

Moreover, \eqref{raimixed} implies $(m_q)^{1/q}\ge A^{-1}m_1\ge A^{-1}>0$ for all $q\ge C$. Since also $(m_q)^{1/q}>0$ for all $1\le q<C$ we get $\liminf_{p\rightarrow\infty}(m_p)^{1/p}>0$, and so \hyperlink{om2}{$(\omega_2)$} for $\omega_M$ (see Lemma \ref{assoweightomega0}).\vspace{6pt}

Consequently, we have $\omega_M\in\hyperlink{omset1}{\mathcal{W}}$ and satisfying \hyperlink{om2}{$(\omega_2)$}. Therefore we can apply \cite[Thm. 3]{characterizationstabilitypaper} to $\omega\equiv\omega_M$ (see also \cite[Thm. 6.3]{compositionpaper}) and \cite[Thm. 5]{characterizationstabilitypaper} to $\mathcal{M}\equiv\mathcal{M}_{\omega_M}$. Finally, by connecting both characterizations via \cite[Thm. 5.14 $(2)$]{compositionpaper}, and for this \hyperlink{om1}{$(\omega_1)$} is indispensable by \eqref{newexpabsorb}, property $(\alpha_0)$ for $\omega_M$ follows and so we are done.\vspace{6pt}

For this note that \hyperlink{om2}{$(\omega_2)$} is assumption $(\omega_4)$ in \cite{compositionpaper} and yields by \cite[Cor. 5.15]{compositionpaper} property $(\mathcal{M}_{\mathcal{H}})$ for $\mathcal{M}_{\omega_M}$ (i.e. $\liminf_{p\rightarrow\infty}(m^{(l)}_p)^{1/p}>0$ for all $l>0$). By \eqref{newmoderategrowth} the matrix $\mathcal{M}_{\omega_M}$ satisfies the mixed derivation closedness properties of both types in \cite{compositionpaper} and \cite{characterizationstabilitypaper}.
\qed\enddemo

\begin{remark}\label{alpha0theoremnewrmark}
Let us mention that in Theorem \ref{alpha0theorem} we can show $(iii)\Rightarrow(ii)$ directly: We consider in \eqref{raimixed} the choices $p\ge 1$ and $q=C\lambda p$ for some $\lambda\in\NN_{>0}$. By multiplying with $p!^{1/p}$ and Stirling's formula we get
\begin{align*}
\lambda(M_p)^{1/p}&=\lambda p!^{1/p}(m_p)^{1/p}\le A \lambda p!^{1/p}(m_{\lambda Cp})^{1/(\lambda C p)}=A\lambda\frac{p!^{1/p}}{(\lambda Cp)!^{1/(\lambda Cp)}}(M_{\lambda Cp})^{1/(\lambda Cp)}
\\&
\le A\lambda e\frac{p}{\lambda Cp}(M_{\lambda Cp})^{1/(\lambda Cp)},
\end{align*}
thus \eqref{raimixedM} is verified for $B:=\frac{Ae}{C}$ (w.l.o.g. we can assume $A\ge C$) and $C=C$.
\end{remark}

When combining Theorem \ref{alpha0theorem} with \cite[Theorem 5.14]{compositionpaper}, \cite[Thm. 6.3, Thm. 6.5]{compositionpaper} and \cite[Thm. 3, Thm. 5]{characterizationstabilitypaper} we get the following information for given $M\in\hyperlink{LCset}{\mathcal{LC}}$:\vspace{6pt}

Some/any of the equivalent conditions in Theorem \ref{alpha0theorem} is satisfied if and only if the class $\mathcal{E}_{[\omega_M]}$ is closed under composition resp. if and only if the Roumieu-type class $\mathcal{E}_{\{\omega_M\}}$ satisfies all further stability properties listed in \cite[Thm. 3, Thm. 5]{characterizationstabilitypaper}. In addition $\mathcal{E}_{[\mathcal{M}_{\omega_M}]}=\mathcal{E}_{[\omega_M]}$ holds as locally convex vector spaces.

However, note that the proof of Theorem \ref{alpha0theorem} is ''tailor-made'' for the Roumieu-type condition \hyperlink{R-rai}{$(\mathcal{M}_{\{\on{rai}\}})$} and that \hyperlink{om5}{$(\omega_5)$}, required as a basic assumption in the Beurling variant result \cite[Thm. 4]{characterizationstabilitypaper}, is not guaranteed in the proof above. In this context we remark that Theorem \ref{alpha0theorem} trivially applies to the constant sequence $m_p=1$ for all $p\in\NN$ which yields for the Roumieu-type the class of real-analytic functions. But it is known that in this case $\omega_M$ is equivalent to $t\mapsto t$ and so clearly violating \hyperlink{om5}{$(\omega_5)$}.

\subsection{Combining $(\alpha_0)$ with moderate growth}\label{constantcase}
The aim is to extend the equivalences listed in Theorem \ref{alpha0theorem} when assuming for $M$ in addition property \hyperlink{mg}{$(\on{mg})$}.

\begin{theorem}\label{alpha0theorem1rem}
Let $M\in\hyperlink{LCset}{\mathcal{LC}}$ be given and satisfying \hyperlink{mg}{$(\on{mg})$}. Then the following properties are equivalent:
\begin{itemize}
\item[$(i)$] $M$ satisfies
\begin{equation}\label{alpha0theorem1equ0mg}
\exists\;B\ge 1\;\forall\;\lambda,p\in\NN_{>0}:\;\;\;(M_p)^{1/p}\lambda\le B(M_{\lambda p})^{1/(\lambda p)}.
\end{equation}
\item[$(ii)$] The sequence of quotients $\mu$ satisfies
\begin{equation}\label{alpha0theorem1equ0mgmu}
\exists\;B\ge 1\;\forall\;\lambda,p\in\NN_{>0}:\;\;\;\mu_p\lambda\le B\mu_{\lambda p}.
\end{equation}

\item[$(iii)$] $M$ is root almost increasing (has \eqref{rootalmost}), so
$$\exists\;A\ge 1\;\forall\;1\le p\le q:\;\;\;(m_p)^{1/p}\le A(m_q)^{1/q}.$$

\item[$(iv)$] $\mu$ satisfies
$$\exists\;A\ge 1\;\forall\;1\le p\le q:\;\;\;\frac{\mu_p}{p}\le A\frac{\mu_q}{q},$$
i.e. the sequence $(\mu_p/p)_p$ is almost increasing.

\item[$(v)$] There exists $S\in\hyperlink{LCset}{\mathcal{LC}}$ such that $M\hyperlink{approx}{\approx}S$ and $s$ is log-convex (i.e. $M$ is equivalent to a strongly log-convex weight sequence $S$).

\item[$(vi)$] $\omega_M$ satisfies $(\alpha_0)$ (\eqref{alpha0}), i.e. is equivalent to a subadditive weight function.

\item[$(vii)$] $\mathcal{M}_{\omega_M}$ satisfies \hyperlink{R-rai}{$(\mathcal{M}_{\{\on{rai}\}})$}.
\end{itemize}
\end{theorem}

Before proving this result we point out:

\begin{itemize}
\item[$(a)$] \eqref{alpha0theorem1equ0mgmu} does imply $(\beta_3)$ (e.g. choose $\lambda=2B$)
and this is consistent with Proposition \ref{beta3comparsion} since assertion $(vi)$ above clearly yields \hyperlink{om1}{$(\omega_1)$} for $\omega_M$.

\item[$(b)$] Assertion $(v)$ is satisfied for any $M\in\hyperlink{LCset}{\mathcal{LC}}$ having \hyperlink{beta1}{$(\beta_1)$} as shown in \cite[Corollary 1.3 $(a)$]{petzsche}. In fact, there the equivalence has been established on the level of the sequence of quotients which is in general stronger than relation \hyperlink{approx}{$\approx$} and both notions do coincide when assuming  \hyperlink{mg}{$(\on{mg})$} (by taking into account \eqref{mgstrange}). But, on the other hand, \hyperlink{mg}{$(\on{mg})$} has not been used in \cite{petzsche}.

\item[$(c)$] The implication $(v)\Rightarrow(vi)$ has been shown in \cite[Lemma 3.4 $(ii)$]{sectorialextensions}: There, on the one hand, one only requires $M\in\RR_{>0}^{\NN}$ but, on the other hand $s\in\hyperlink{LCset}{\mathcal{LC}}$ and so in addition $(s_k)^{1/k}\rightarrow+\infty$ as $k\rightarrow+\infty$. But this, by equivalence, would yield $(m_k)^{1/k}\rightarrow+\infty$ as well.
\end{itemize}

\demo{Proof}
First, we are showing that \eqref{raimixedM} is equivalent to \eqref{alpha0theorem1equ0mg}.

The implication \eqref{alpha0theorem1equ0mg}$\Rightarrow$\eqref{raimixedM} is clear since $p\mapsto(M_p)^{1/p}$ is nondecreasing. For the converse, we require that
$$\forall\;C\in\NN_{>0}\;\exists\;A\ge 1\;\forall\;\lambda,p\in\NN_{>0}:\;\;\;(M_{\lambda Cp})^{1/(\lambda Cp)}\le A(M_{\lambda p})^{1/(\lambda p)},$$
which is equivalent to $M_{\lambda Cp}\le A^{\lambda Cp}(M_{\lambda p})^{C}$. This follows directly by iterating property \hyperlink{mg}{$(\on{mg})$} (note that the number of iterations is only depending on the constant $C$).\vspace{6pt}

$(i)\Leftrightarrow(vi)$ We apply Theorem \ref{alpha0theorem} and the shown equivalence just before.\vspace{6pt}

$(vi)\Leftrightarrow(vii)$ This is shown in Theorem \ref{alpha0theorem}.\vspace{6pt}

$(i)\Leftrightarrow(ii)$ This is immediate by log-convexity, normalization and \hyperlink{mg}{$(\on{mg})$} by involving \eqref{mgstrange}.\vspace{6pt}

$(i)\Leftrightarrow(iii)$ This follows by combining \eqref{raimixedM}$\Leftrightarrow$\eqref{alpha0theorem1equ0mg}, $(ii)\Leftrightarrow(iv)$ in Theorem \ref{alpha0theorem} and $(ii)\Leftrightarrow(iii)$ in Corollary \ref{lccorollary}.\vspace{6pt}

$(iii)\Leftrightarrow(iv)$ By log-convexity, normalization and \hyperlink{mg}{$(\on{mg})$} we recall that (see \eqref{mgstrangem})
$$\exists\;A\ge 1\;\forall\;p\ge 1:\;\;\;(m_p)^{1/p}\le e\frac{\mu_p}{p}\le eA(m_p)^{1/p}.$$
This implies the desired equivalence.\vspace{6pt}

$(iv)\Rightarrow(v)$ This has been shown in \cite[Lemma 8]{whitneyextensionmixedweightfunctionII} (inspired by \cite[Prop. 4.15]{JimenezGarridoSanz}); more precisely we define $S$ in term of the quotient sequence $\sigma$ by taking $\sigma_p:=p\inf_{q\ge p}\frac{\mu_q}{q}$, $p\ge 1$, and $\sigma_0:=1$.\vspace{6pt}

$(v)\Rightarrow(iv)$ By $M\hyperlink{approx}{\approx}S$ we get $\frac{1}{C}(m_p)^{1/p}\le(s_p)^{1/p}\le C(m_p)^{1/p}$ for some $C\ge 1$ and all $p\in\NN_{>0}$. Moreover, since $M$ does also have \hyperlink{mg}{$(\on{mg})$}, by \eqref{mgstrangem} we get
$$\exists\;A,C\ge 1\;\forall\;1\le p\le q:\;\;\;\frac{\mu_p}{p}\le A(m_p)^{1/p}\le AC(s_p)^{1/p}\le AC(s_q)^{1/q}\le AC^2(m_q)^{1/q}\le AC^2e\frac{\mu_q}{q}.$$
\qed\enddemo

{\itshape Note:} In the situation of Theorem \ref{alpha0theorem1rem}, if some/any of the equivalent conditions there is satisfied, then additionally to the comments on the stability properties for $\mathcal{E}_{[\omega_M]}$ and $\mathcal{E}_{\{\omega_M\}}$ mentioned at the end of the previous section, we get that $\mathcal{E}_{[M]}=\mathcal{E}_{[\mathcal{M}_{\omega_M}]}=\mathcal{E}_{[\omega_M]}$ as locally convex vector spaces (by applying the results from \cite{BonetMeiseMelikhov07} and \cite[Cor. 5.8 $(2)$, Thm. 5.14 $(3)$]{compositionpaper}).



\subsection{Stability under equivalence relations}\label{stabilitysection}
In this section we comment on the stability of $(\alpha_0)$ w.r.t. the relevant equivalence relations \hyperlink{approx}{$\approx$} and \hyperlink{sim}{$\sim$}. Recall that \hyperlink{approx}{$\approx$} is characterizing the equivalence of the corresponding ultradifferentiable function classes $\mathcal{E}_{[M]}$, whereas \hyperlink{sim}{$\sim$} becomes relevant for $\mathcal{E}_{[\omega]}$, see \cite[Prop. 2.12 $(1)$, Cor. 5.17 $(1)$]{compositionpaper}.\vspace{6pt}

Clearly, $(\alpha_0)$ is preserved under \hyperlink{sim}{$\sim$} and the aim is to show that it is also preserved under  \hyperlink{approx}{$\approx$} when considered for associated weight functions.

\begin{corollary}\label{alpha0theorem1remcor}
Let $M\in\hyperlink{LCset}{\mathcal{LC}}$ be given such that $M\hyperlink{approx}{\approx}N$. Assume that $\omega_M$ has $(\alpha_0)$, then we get for any $N\in\hyperlink{LCset}{\mathcal{LC}}$ satisfying $M\hyperlink{approx}{\approx}N$ that

\begin{itemize}
\item[$(i)$] $\omega_N$ has $(\alpha_0)$ and

\item[$(ii)$] $\omega_M\hyperlink{sim}{\sim}\omega_N$ holds true.
\end{itemize}
In particular, in the situation of Theorem \ref{alpha0theorem1rem} we get that $\omega_S\hyperlink{sim}{\sim}\omega_M$ and $\omega_S$ has $(\alpha_0)$, too.
\end{corollary}

\demo{Proof}
$(i)$ Property \eqref{raimixedM} is clearly preserved under \hyperlink{approx}{$\approx$} by changing the constant $B$ (resp. \eqref{raimixed} is preserved by changing $A$). So Theorem \ref{alpha0theorem} yields the conclusion.\vspace{6pt}

$(ii)$ By $(iii)$ in Theorem \ref{alpha0theorem} we get \hyperlink{om1}{$(\omega_1)$} for $\omega_M$. This fact and $M\hyperlink{approx}{\approx}N$ yield $\omega_M\hyperlink{sim}{\sim}\omega_N$, see \cite[Remark 3.3]{sectorialextensions}.
\qed\enddemo

\section{Comments on condition $(\alpha_1)$}\label{alpha1section}
In this section we give some comments on $(\alpha_1)$. This is for the sake of completeness, in order to illustrate the subtle difference between conditions $(\alpha_0)$ and $(\alpha_1)$ for associated weight functions and to comment on the differences in the proofs. First note that $(\alpha_1)$ is equivalent to having
\begin{equation}\label{alpha1rewritten}
\exists\;C\ge 1\;\forall\;\lambda\ge 1\;\exists\;D_{\lambda}\ge 1\;\forall\;t\ge 0:\;\;\;\omega(\lambda t)\le C\lambda\omega(t)+D_{\lambda}\lambda.
\end{equation}
In \eqref{alpha1rewritten} equivalently we could replace $D_{\lambda}\lambda$ by $D'_{\lambda}$ (as it has been done in \cite{BonetDomanski07} and \cite{BonetDomanski06}) but in analogy to \eqref{alpha0rewritten} we prefer to write this condition in this way. Also $(\alpha_1)$ implies \hyperlink{om1}{$(\omega_1)$} by choosing $\lambda=2$.\vspace{6pt}

Let us now formulate the following characterization of $(\alpha_1)$ analogous to Theorem \ref{alpha0theorem}, here we can even treat a mixed setting.

\begin{theorem}\label{alpha0theoremcor}
Let $M,N\in\hyperlink{LCset}{\mathcal{LC}}$ be given and let $\omega_M,\omega_N$ be the corresponding associated weight functions. Then the following are equivalent:

\begin{itemize}
\item[$(i)$] The associated weight functions satisfy
\begin{equation}\label{alpha1mixed}
\exists\;A\ge 1\;\exists\;B\ge 1\;\exists\;C\ge 1\;\forall\;\lambda\ge 1\;\exists\;D\ge 1\;\forall\;t\ge 0:\;\;\;\omega_N(\lambda Ct)\le B\lambda\omega_M(At)+D\lambda,
\end{equation}
i.e. the mixed version of $(\alpha_1)$ (\eqref{alpha1rewritten}).

\item[$(ii)$] $M$ and $N$ satisfy
\begin{equation}\label{omega1mixed}
\exists\;C\in\NN_{>0}\;\exists\;B\ge 1\;\forall\;\lambda\in\NN_{>0}\;\exists\;D\ge 1\;\forall\;p\ge 1:\;\;\;(M_p)^{1/p}\lambda\le BD^{1/p}(N_{\lambda Cp})^{1/(\lambda Cp)}.
\end{equation}

\item[$(iii)$] $m$ and $n$ satisfy
\begin{equation*}\label{omega1mixedmnew}
\exists\;C\in\NN_{>0}\;\exists\;A\ge 1\;\forall\;\lambda\in\NN_{>0}\;\exists\;D\ge 1\;\forall\;p\ge 1:\;\;\;(m_p)^{1/p}\le AD^{1/p}(n_{\lambda Cp})^{1/(\lambda Cp)}.
\end{equation*}
\end{itemize}
\end{theorem}

\demo{Proof}
$(i)\Rightarrow(ii)$ and $(ii)\Rightarrow(iii)$ follow analogously as in Theorem \ref{alpha0theorem}, $(iii)\Rightarrow(ii)$ as in Remark \ref{alpha0theoremnewrmark}. (Note that all these implications are also valid in the mixed setting.)\vspace{6pt}

$(ii)\Rightarrow(i)$ We use the same trick as for $(i)\Rightarrow(ii)$ in Theorem \ref{omega1charact}. By using \eqref{Mxonestable} and \eqref{Mgoodtransform} we get $(N_{\lambda Cp})^{1/(\lambda C)}=(N^{(1)}_{\lambda Cp})^{1/(\lambda C)}=N^{(\lambda C)}_p$, with $\mathcal{M}_{\omega_N}:=\{N^{(l)}: l>0\}$ denoting the matrix associated with the weight $\omega_N$. Then \eqref{omega1mixed} implies (after multiplying with $t\ge 0$) that
$$\exists\;C\in\NN_{>0}\;\exists\;B\ge 1\;\forall\;\lambda\in\NN_{>0}\;\exists\;D_{\lambda}\ge 1\;\forall\;p\in\NN\;\forall\;t\ge 0:\;\;\;\frac{(\lambda t)^p}{N^{(\lambda C)}_p}\le D_{\lambda}\frac{(Bt)^p}{M_p},$$
yielding $\omega_{N^{(\lambda C)}}(\lambda t)\le\omega_M(Bt)+D'_{\lambda}$ by definition of associated weight functions and we have set $D'_{\lambda}:=\log(D_{\lambda})$. Thus \eqref{goodequivalence} implies
\begin{align*}
&\exists\;C\in\NN_{>0}\;\;\exists\;B\ge 1\;\forall\;\lambda\in\NN_{>0}\;\exists\;D'_{\lambda}\ge 1\;\exists\;D_{\lambda C}\ge 1\;\forall\;t\ge 0:
\\&
\frac{1}{2\lambda C}\omega_{N}(\lambda t)-\frac{D_{\lambda C}}{2\lambda C}=\frac{1}{2\lambda C}\omega_{N^{(1)}}(\lambda t)-\frac{D_{\lambda C}}{2\lambda C}\le\omega_{N^{(\lambda C)}}(\lambda t)\le\omega_M(Bt)+D'_{\lambda},
\end{align*}
hence
\begin{align*}
&\exists\;C\in\NN_{>0}\;\exists\;B\ge 1\;\forall\;\lambda\in\NN_{>0}\;\exists\;D'_{\lambda}\ge 1\;\exists\;D_{\lambda C}\ge 1\;\forall\;t\ge 0:
\\&
\omega_N(\lambda t)\le 2C\lambda\omega_M(Bt)+2C D'_{\lambda}\lambda+D_{\lambda C}.
\end{align*}
This proves \eqref{alpha1mixed} for all $\lambda\in\NN_{>0}$ (by taking $C=1$, $B=2C$, $A=B$, $D=\max\{2CD'_{\lambda},D_{\lambda C}\}$).

In order to finish, let now $\lambda\ge 1$ be arbitrary, then $\mu\le\lambda<\mu+1$ for some $\mu\in\NN_{>0}$. Hence
\begin{align*}
\omega_N(\lambda t)&\le\omega_N((\mu+1)t)\le 2C(\mu+1)\omega_M(Bt)+2C D'_{\mu+1}(\mu+1)+D_{(\mu+1)C}
\\&
\le 4C\mu\omega_M(Bt)+4C D'_{\mu+1}\mu+D_{(\mu+1)C}\le 4C\lambda\omega_M(Bt)+4CD'_{\mu+1}\lambda+D_{(\mu+1)C},
\end{align*}
and it suffices to take $C=1$, $B=4C$, $A=B$, $D=\max\{4CD'_{\lceil\lambda\rceil}, D_{\lceil\lambda\rceil C}\}$.
\qed\enddemo

When taking $M=N$ in Theorem \ref{alpha0theoremcor}, then we obtain the following characterization:

\begin{corollary}\label{alpha0theoremcorcor}
Let $M\in\hyperlink{LCset}{\mathcal{LC}}$ be given. Then $\omega_M$ satisfies $(\alpha_1)$ (\eqref{alpha1rewritten}) if and only if some/any of the equivalent conditions in Theorem \ref{alpha0theoremcor} holds with $M=N$.
\end{corollary}

\demo{Proof}
It remains to show that \eqref{alpha1rewritten} for $\omega_M$ is equivalent to \eqref{alpha1mixed}.

First, \eqref{alpha1rewritten} implies \eqref{alpha1mixed} for $M=N$ with $A=C=1$ and $B=C$.

For the converse, first note that some/any of the equivalent conditions in Theorem \ref{alpha0theoremcor} implies \hyperlink{om1}{$(\omega_1)$} for $\omega_M$: In \eqref{omega1mixed} take $\lambda\in\NN_{>0}$ with $\lambda>B$ ($B$ is not depending on $\lambda$) and get $\frac{(M_{\lambda Cp})^{1/(\lambda Cp)}}{(M_p)^{1/p}}\ge\frac{\lambda}{B}\frac{1}{D^{1/p}}\rightarrow\frac{\lambda}{B}>1$ as $p\rightarrow\infty$, hence \eqref{omega1charactequalternative} is verified with $L:=\lambda C$.

So we set $s:=Ct$ in \eqref{alpha1mixed} which gives $\omega_M(\lambda s)\le B\lambda\omega_M(sA/C)+D\lambda$. When $A\le C$, then we are done since $\omega_M$ is nondecreasing. If $C<A$, then take $n\in\NN$ minimal such that $2^n\ge\frac{A}{C}$ and then by iterating \hyperlink{om1}{$(\omega_1)$} we have
$$\exists\;L'\ge 1\;\forall\;s\ge 0:\;\;\;\omega_M(\lambda s)\le B\lambda\omega_M(2^ns)+D\lambda\le B L'\lambda\omega_M(s)+BL'\lambda+D\lambda.$$
Note that both $A$ and $C$ are not depending on $\lambda$, hence $n$ and so the new constant $L'$ will be not depending as well.
\qed\enddemo

When comparing the proofs of Theorem \ref{alpha0theoremcor} and Theorem \ref{alpha0theorem} we can see the following differences concerning the used techniques:

\begin{itemize}
\item[$(a)$] Inspecting the proof of $(ii)\Rightarrow(i)$ in Theorem \ref{alpha0theoremcor} we see that even when assuming the stronger assumption \eqref{raimixedM} we cannot infer \eqref{alpha0mixed}, again only the weaker condition \eqref{alpha1mixed}: More precisely, we still have the dependence of $\lambda$ by the constant $D_{\lambda C}$ (by involving \eqref{goodequivalence}). By inspecting the proof given in \cite[Lemma 2.5]{sectorialextensions} (and in \cite[Thm. 4.0.3]{dissertation}) we see that $D_{\lambda C}=\lambda Cd_{\lambda C}$ for some different constant $d_{\lambda C}$ which may, in general, tend to infinity as $\lambda\rightarrow\infty$.

\item[$(b)$] In the proof of Theorem \ref{alpha0theorem} we have been able to overcome this problem by involving the characterizations from \cite{characterizationstabilitypaper} and \cite{compositionpaper}. However, in the mixed setting this idea failed and note that the proof of Theorem \ref{alpha0theoremcor} does not require any additional information on the underlying ultradifferentiable function classes.
\end{itemize}

As already mentioned in the introduction the weight
\begin{equation}\label{omegacounter}
\widetilde{\omega}(t):=0,\;\;\;0\le t\le 1,\hspace{30pt}\widetilde{\omega}(t):=t\log(t),\;\;\;t\ge 1,
\end{equation}
has $(\alpha_1)$ but not $(\alpha_0)$. Note that $\widetilde{\omega}\in\hyperlink{omset1}{\mathcal{W}}$ and \hyperlink{om6}{$(\omega_6)$} is valid, however the standard assumption \hyperlink{om2}{$(\omega_2)$}, and hence \hyperlink{om5}{$(\omega_5)$}, is violated. So, from this point of view, the difference between $(\alpha_0)$ and $(\alpha_1)$ can be considered as large and in general $\omega_M$ satisfying \eqref{alpha1mixed} (with $M=N$) will violate \hyperlink{om2}{$(\omega_2)$}.

Similar arguments as given for $(\alpha_0)$ in Section \ref{stabilitysection} are valid for $(\alpha_1)$ as well.

\vspace{12pt}
{\itshape Open problem:}

Construct $M\in\hyperlink{LCset}{\mathcal{LC}}$ such that $\omega_M$ has $(\alpha_1)$ but not $(\alpha_0)$; thus in view of Theorems \ref{alpha0theoremcor} and \ref{alpha0theorem} that $M$ has \eqref{omega1mixed} (with $M=N$) but not \eqref{raimixedM}.

By Lemma \ref{assoweightomega0} we are also interested in having $\liminf_{p\rightarrow\infty}(m_p)^{1/p}>0$ which implies \hyperlink{om2}{$(\omega_2)$} for $\omega_M$ since $\widetilde{\omega}$ from \eqref{omegacounter} does violate this basic assumption.

\section{Strong sufficient conditions for having subadditivity-like assumptions}\label{addcommsection}
Recall that for the conditions studied so far we have the following implications:
$$(\alpha_0)\Longrightarrow(\alpha_1)\Longrightarrow\hyperlink{om1}{(\omega_1)}.$$
For reasons of completeness, we give now some comments on conditions which ''prolong'' this chain on the left-hand side.\vspace{6pt}

First, we recall the definition of $\gamma(\omega)$, $\omega$ a given weight function (see \cite[Sect. 2.3]{index} and the references therein): Let $\gamma>0$, then we say that $\omega$ has property $(P_{\omega,\gamma})$ if
\begin{equation*}\label{newindex1}
\exists\;K>1:\;\;\;\limsup_{t\rightarrow+\infty}\frac{\omega(K^{\gamma}t)}{\omega(t)}<K.
\end{equation*}
If $(P_{\omega,\gamma})$ holds for some $K>1$, then also $(P_{\omega,\gamma'})$ is satisfied for all $\gamma'\le\gamma$ with the same $K$. Moreover we can restrict to $\gamma>0$, because for $\gamma\le 0$ condition $(P_{\omega,\gamma})$ is satisfied for all weights $\omega$ (since $\omega$ is nondecreasing and $K>1$). Then we put
\begin{equation}\label{newindex2}
\gamma(\omega):=\sup\{\gamma>0: (P_{\omega,\gamma})\;\;\text{is satisfied}\}.
\end{equation}

First, we recall that $\gamma(\omega)>0$ if and only if \hyperlink{om1}{$(\omega_1)$} holds true, see \cite[Cor. 2.14]{index}. By \cite[Thm. 2.11, Cor. 2.13]{index} we have that $\gamma(\omega)>1$ if and only if $\omega$ has
\begin{equation}\label{assostrongnq}
\exists\;C\ge 1\;\forall\;y\ge 0:\;\;\;\int_1^{+\infty}\frac{\omega(yt)}{t^2}dt\le C\omega(y)+C,
\end{equation}
i.e. $\omega$ has the {\itshape strong non-quasianalyticity} condition for weight functions. It is known that each weight function satisfying \eqref{assostrongnq} is equivalent to a concave weight, more precisely to $\kappa_{\omega}(y):=\int_1^{+\infty}\frac{\omega(yt)}{t^2}dt$, see \cite[Prop. 1.3]{MeiseTaylor88}. Summarizing, we get that $\gamma(\omega)>1$ implies $(\alpha_0)$.\vspace{6pt}

Next we introduce the following growth assumption on $\omega$:
\begin{equation}\label{omega7}
\exists\;H>0:\;\;\;\omega(t^2)=O(\omega(Ht)),\;\;\;\text{as}\;\;t\rightarrow+\infty.
\end{equation}
We refer to \cite[Appendix A]{sectorialextensions1}, where this condition is denoted by $(\omega_7)$, and in \cite[Lemma A.1]{sectorialextensions1} we have seen that this condition implies $\gamma(\omega)=+\infty$ and mentioned that this implication is in general strict.

Gathering everything, we have the following chain of implications for weight functions:
\begin{equation}\label{chain}
\eqref{omega7}\Longrightarrow\gamma(\omega)=+\infty\Longrightarrow\gamma(\omega)>1\Longrightarrow(\alpha_0)\Longrightarrow(\alpha_1)\Longrightarrow\hyperlink{om1}{(\omega_1)}\Longleftrightarrow\gamma(\omega)>0.
\end{equation}

We focus now on $\omega\equiv\omega_M$, $M\in\hyperlink{LCset}{\mathcal{LC}}$. By combining \cite[Cor. 2.14]{index} with Theorem \ref{omega1charact} we immediately get
$$\gamma(\omega_M)>0\Longleftrightarrow\omega_M\;\text{has}\;\hyperlink{om1}{(\omega_1)}\Longleftrightarrow\exists\;L\in\NN_{>0}:\;\;\;\liminf_{p\rightarrow+\infty}\frac{(M_{Lp})^{\frac{1}{Lp}}}{(M_p)^{\frac{1}{p}}}>1.$$
By \cite[Prop. 4.4]{Komatsu73} we have that \hyperlink{beta1}{$(\beta_1)$} implies \eqref{assostrongnq} for $\omega_M$ but in general the converse is not true and the difference between both notions of strong non-quasianalyticity can be large, see \cite[Sect. 5]{index}.\vspace{6pt}

Concerning \eqref{omega7} for $\omega_M$, in \cite[Corollary 6.4]{scales} by using similar techniques as in this paper, the following result is shown:

\begin{proposition}\label{omega7prop}
Let $M\in\hyperlink{LCset}{\mathcal{LC}}$ be given. Then $\omega_M$ does have \eqref{omega7} if and only if
\begin{equation}\label{omega7propequ}
\exists\;L\in\NN_{>0}\;\exists\;A,B\ge 1\;\forall\;p\in\NN:\;\;\;(M_p)^{2L}\le AB^pM_{Lp}.
\end{equation}
In fact, we have to choose $L\ge 2$, otherwise this would yield a contradiction to $\lim_{p\rightarrow+\infty}(M_p)^{1/p}=+\infty$.
\end{proposition}

Then \eqref{omega7propequ} yields $\frac{(M_{Lp})^{1/(Lp)}}{(M_p)^{1/p}}\ge\frac{(M_p)^{1/p}}{A^{1/(Lp)}B^{1/L}}$ and so one has
$$\exists\;L\in\NN_{>0}:\;\;\;\lim_{p\rightarrow+\infty}\frac{(M_{Lp})^{1/(Lp)}}{(M_p)^{1/p}}=+\infty.$$
Consequently, \eqref{omega7propequ} does imply \eqref{omega1charactequalternative} and this is consistent with Theorem \ref{omega1charact} and with \eqref{chain} applied to $\omega_M$. Thus, also in the weight sequence setting, one can view \eqref{omega1charactequalternative} as the ''weakest'' and \eqref{omega7propequ} as the ''strongest'' growth condition when dealing with expressions of the form $\frac{(M_{Lp})^{1/(Lp)}}{(M_p)^{1/p}}$.\vspace{6pt}

Finally, let us mention that the conclusion of Proposition \ref{beta3comparsion} also follows by combining \cite[Cor. 2.14, Thm. 3.11 $(v)\Leftrightarrow(vii)$, Cor. 4.6 $(i),(iii)$]{index}. The arguments in \cite{index} are completely different and relying on the notion of growth indices. For this also the index $\gamma(M)$ is showing up, see \cite[Sect. 1.3]{Thilliezdivision} and \cite[Sect. 3.1]{index}, and it is compared with $\gamma(\omega_M)$. Both indices $\gamma(\omega)$ and $\gamma(M)$ are needed for extension results in the ultraholomorphic setting.

\bibliographystyle{plain}
\bibliography{Bibliography}

\end{document}